\documentclass[10pt]{article}

\usepackage{mathrsfs}
\usepackage{amsfonts}
\usepackage{amssymb}

\font\pass cmss10  at 6.9pt
\font\passauthor cmr10  at 8pt
\font\passtitle cmti10  at 8pt

\newenvironment{proof}{\noindent {\it Proof.}}{\hfill$\Box$\bigskip}

\newenvironment{passage}
  {\list{}{\listparindent 0.2cm%
  \baselineskip 0pt
  \setlength{\leftmargin}{4.5cm}
  \setlength{\rightmargin}{-0.5cm}
  }%
  \item\relax}
  {\endlist}

\DeclareSymbolFont{EUR}{U}{eur}{m}{n}
\SetSymbolFont{EUR}{bold}{U}{eur}{b}{n}
\DeclareSymbolFontAlphabet{\eur}{EUR}

\DeclareSymbolFont{EUB}{U}{eur}{b}{n}
\SetSymbolFont{EUB}{bold}{U}{eur}{b}{n}
\DeclareSymbolFontAlphabet{\eub}{EUB}

\DeclareSymbolFont{EUS}{U}{eus}{m}{n}
\SetSymbolFont{EUS}{bold}{U}{eus}{b}{n}
\DeclareSymbolFontAlphabet{\eus}{EUS}

\DeclareSymbolFont{EUF}{U}{euf}{m}{n}
\SetSymbolFont{EUF}{bold}{U}{euf}{b}{n}
\DeclareSymbolFontAlphabet{\euf}{EUF}

\DeclareSymbolFont{AMSb}{U}{msb}{m}{n}
\DeclareSymbolFontAlphabet{\mathbb}{AMSb}

\newcommand\rank{\mathop{\rm rank}}

\newcommand\p{\partial}
\newcommand{\at}[1]{\vert\sb{\sb{#1}}}

\def\R{\mathbb{R}}
\def\S{\mathbb{S}}
\def\hvar{{\hbar}}
\newcommand\C{{\mathbb C}}
\newcommand\N{{\mathbb N}}\newcommand\Z{{\mathbb Z}}
\newcommand{\Abs}[1]{\left\vert#1\right\vert}
\newcommand{\abs}[1]{\vert #1 \vert}

\newcommand{\norm}[1]{\Vert #1 \Vert}

\newcommand\const{\mathop{\rm const}}

\newcommand\sothat{\,\,{\rm :}\,\,}

\def\hardy{\euf{h}}
\def\bmo{{\rm BMO}}

\def\varDelta{{\mathit\Delta}}

\def\varSigma{\mathit\Sigma}
\def\text#1{{\rm #1}}

\newtheorem{theorem}{Theorem}[section]

\newtheorem{lemma}[theorem]{Lemma}
\newtheorem{corollary}[theorem]{Corollary}
\newtheorem{proposition}[theorem]{Proposition}

\newtheorem{definition}[theorem]{Definition}
\newtheorem{remark}[theorem]{Remark}
\newtheorem{example}[theorem]{Example}

\newcommand{\sect}[1]{\setcounter{equation}{0}\section{#1}}

\begin{document}


\title{$L\sp p\mbox{-}L\sp q$ regularity of Fourier integral operators with caustics}

\author{Andrew Comech
\footnote{This 
work was supported in part by the NSF under Grants No. 0296036 and 0200880}
}

\maketitle


\begin{abstract}
{\small
The caustics of Fourier integral operators
are defined as caustics of the corresponding Schwartz kernels 
(Lagrangian distributions on $X\times Y$).
The caustic set $\Sigma(\eub C)$
of the canonical relation $\eub C$
is characterized as the set of points
where the rank of the projection
$\pi:\eub C\to X\times Y$
is smaller than its maximal value, $\dim(X\times Y)-1$.
We derive the $L\sp p(Y)\to L\sp q(X)$ estimates
on Fourier integral operators with caustics of corank $1$
(such as caustics of type $A\sb{m+1}$, $m\in\N$).
For the values of $p$ and $q$
outside of certain neighborhood of 
the line of duality, $q=p'$,
the $L\sp p\to L\sp q$ estimates
are proved to be caustics-insensitive.

\noindent
We apply our results to the analysis 
of the blow-up of the estimates
on the half-wave operator
just before the geodesic flow forms caustics.
}
\end{abstract}

\begin{passage}

\noindent
\pass
Stretching his stiff legs, Persikov got up, returned to his laboratory,
yawned, rubbed his permanently inflamed eyelids, sat down on the stool and
looked into the microscope .\,.\,.\,\,With his right eye Persikov saw the cloudy white plate and
blurred pale amoebas on it, but in the middle of the plate sat a coloured
tendril, like a female curl. 
The coloured streak of light merely got in the way and indicated
that the specimen was out of focus. 
The zoologist's long fingers had already tightened on the
knob, when suddenly they trembled and let go .\,.\,.\,\,He noticed 
that one particular ray in the coloured tendril 
stood out more vividly and boldly than the others .\,.\,.\,\,This strip of red was teeming with life. 
The old amoebas were forming pseudopodia in a desperate effort to reach the red strip, 
and when they did they came to life, as if by magic. 
They split into two in the ray, and each of the parts became a new, fresh organism in a 
couple of seconds. 

\medskip

\hfill \passauthor Mikhail Bulgakov, {\passtitle The Fateful Eggs}, 1924.

\end{passage}
				   
\sect{Introduction}
Caustics are the envelopes of the light rays.
At the caustic points, intensity of light is singularly large,
causing different physical phenomena (such as the one 
observed by Professor Persikov).
Mathematically, caustics could be characterized as
points where usual bounds on oscillatory
integrals are no longer valid.
In this paper we are going to consider
how this concept applies to
Fourier integral operators.
This question becomes interesting
in view of a recent paper \cite{MR2000i:35115}
on dissipative semilinear oscillations,
where the $L\sp q$ estimates
on oscillatory integrals with caustics
played the central role.
Our goal is to investigate
how the regularity properties of Fourier integral operators
are affected by the presence of caustics.
We will show, in particular,
that for $q$ away from a certain neighborhood of $q=p'$
the $L\sp p\to L\sp q$ bounds
on Fourier integral operators
are caustic-insensitive.

Oscillatory integrals with caustics have enjoyed much attention.
The classical references are \cite{MR89g:58024}
and \cite{MR53:9306}.
The asymptotics of oscillatory integrals
near caustics were derived in 
\cite{MR33:4446} and \cite{MR58:24404}.

Let us mention previously known estimates
on Fourier integral operators.
The $L\sp 2$ estimates on Fourier integral operators
were considered by H\"ormander \cite{MR52:9299}.
The $L\sp p\to L\sp q$
and $L\sp p\to L\sp{p'}$
estimates on Fourier integral operators
in the context of strictly hyperbolic 
equations with constant coefficients
were addressed in
\cite{MR41:876},
\cite{MR50:10909},
\cite{MR52:8658},
and
\cite{MR95j:35128, MR97f:35122, MR2000i:35113}.
$L\sp p\to L\sp {p'}$ estimates
for certain hyperbolic equations with smooth coefficients
and applications to 
the existence and uniqueness results
for semi-linear hyperbolic equations
are in \cite{MR55:3877}.
The $L\sp p\to L\sp p$ estimates
were derived by Seeger, Sogge, and Stein \cite{MR92g:35252}.
For more information on regularity properties
of generalized Radon transforms
and Fourier integral operators
associated to local graphs and
to degenerate canonical relations
see the reviews \cite{MR2001j:58047, MR1964817}.

\bigskip

We first recall some background 
about caustics of oscillatory integrals.
Let us consider an oscillatory integral
\begin{equation}
u\sb\tau(x)=\tau\sp{k/2}\int\sb{\R^k}
e\sp{i\tau\psi(x,\alpha)}a(x,\tau,\alpha)\,d^k \alpha,
\qquad x\in\R^n,\quad \alpha\in\R^k,\quad\tau>0.
\end{equation}
We assume that
$\psi$ is a smooth function
and that $a\in S^d$ is a symbol
of order $d$ in $\tau$,
compactly supported in $\alpha$ and $x$.
If there are no critical points
of the map $\alpha\mapsto\psi(x,\alpha)$,
so that $\psi\sb\alpha'\ne 0$
everywhere in an open neighborhood
of the support of $a(x,\tau,\alpha)$,
then the repeated integration by parts shows that
$\abs{u\sb\tau(x)}=o(\tau\sp{-N})$, for any $N>0$.
If there are non-degenerate critical points,
where $\psi\sb\alpha'=0$
but $\det\psi\sb{\alpha\alpha}''\ne 0$,
then the method of stationary phase shows that
$\abs{u\sb\tau(x)}=O(\tau\sp d)$.
If we also assume that
$\rank \psi\sb{x\alpha}''\ge n$
(when $k\ge n$),
then one can readily show that
$\norm{u\sb\tau(x)}\sb{L\sp 2}=O(\tau\sp d)$.
It follows that as long as all critical points
are non-degenerate, $u\sb\tau(x)\in L\sp q(\R^n)$, $2\le q\le\infty$,
with the norms bounded uniformly in $\tau\in(0,\infty)$.

If there are degenerate critical points, known as caustics,
then $\norm{u\sb\tau(x)}\sb{L\sp\infty}$
is no longer uniformly bounded.
The order of a caustic $\kappa$ is defined
as the infimum of $\kappa'$ so that
$\norm{u\sb\tau(x)}\sb{L\sp\infty}=O(\tau\sp{\kappa'})$.
For example,
$\psi(x,\alpha)=\alpha^3+x\alpha$
corresponds to the fold ($A\sb 2$),
with $\kappa=1/6$;
$\psi(x,\alpha)=\alpha^4+x\sb 1\alpha^2+x\sb 2$
corresponds to the cusp ($A\sb 3$),
with $\kappa=1/4$.
For more details,
see \cite{MR89g:58024}, \cite{MR53:9306}, \cite{MR96m:58245}.
At the same time,
it was shown in \cite{MR2000i:35115} that
there exists $q\sb c>2$ such that the $L\sp q$ estimates
for $2\le q< q\sb c$ are still bounded
uniformly in $\tau$.
(This information was used to deduce that
the singularities of solutions to
dissipative semilinear equations, like
$\Box u=\dot u\abs{\dot u}\sp{p-1}$,
are absorbed at the caustic if $p$ is larger
than certain critical value;
for generic caustics, one needs $p\ge 3$
for such an absorption to take place.)

\bigskip

Now we turn to Fourier integral operators.
Let $X$ and $Y$ be two smooth manifolds (without boundary).
A Fourier integral operator
$\euf{F}:C\sp\infty\sb{0}(Y)\to \mathscr{D}'(X)$
can be defined (locally) by
\begin{equation}\label{generic-fio}
\euf{F} u(x)=\int\sb{\R^N\times Y}
e\sp{i\phi(x,\theta,y)}a(x,\theta,y)u(y)\,d\theta\,dy,
\end{equation}
where
$a$ is a symbol of order $d$
and $\phi$ is a non-degenerate phase function.
We write $\euf{F}\in I\sp\mu(X,Y,\eub C)$,
where the order of the operator
is defined by $\mu=d+\frac{N}{2}-\frac{\dim X+\dim Y}{4}$
and $\eub C$ is the associated canonical relation.
We will always assume that
\[
\dim X=\dim Y=n
\]
and that
the symbol $a$ is compactly supported in $X\times Y$.

Let us consider $L\sp 1\to L\sp\infty$ estimates
on $\euf{F}$. From (\ref{generic-fio})
one can see that
\begin{equation}\label{1-infinity-condition}
\euf{F}:L\sp 1(Y)\to L\sp\infty(X)
\end{equation}
if $d+N<0$ (which is equivalent with $\mu<-(n+N)/2$).
The smaller the minimal number of oscillatory variables
is, the better $L\sp 1\to L\sp\infty$ regularity properties
$\euf{F}$ possesses.
As we know
from \cite{MR52:9299},
the minimal number of oscillatory variables is equal to
$N\sb{\rm min}=2n-r$, where $r$
is the minimal value of the rank of the projection
$\pi\sb{X\times Y}$ from $\eub C$ onto $X\times Y$:
\[
r=\min\rank d\pi\sb{X\times Y}.
\]
We define the caustic set of the canonical relation
as a subset $\Sigma(\eub C)$ of $\eub C$ where
the rank of $d\pi\sb{X\times Y}$ is not maximal:
\[
\Sigma(\eub C)
=\{
p\in\eub C
\sothat\rank d\pi\sb{X\times Y}\at{\Sigma(\eub C)}<2n-1
\},
\]
so that outside of $\Sigma(\eub C)$
the number of oscillatory variables
of a Fourier integral operator
$\euf{F}\in I\sp\mu(X,Y,\eub C)$
could be reduced to $N=1$.
Let $\euf{F}\in I^\mu(X,Y,\eub C)$,
and let $\euf{F}\sb\lambda$, $\lambda\ge 1$,
be its Littlewood-Paley decomposition.
Similarly to \cite{MR96m:58245},
we will say that $\kappa$ is the highest order
of caustics of $\eub C$
if it is the infimum of numbers $\kappa'$
such that the Schwartz kernel of $\euf{F}\sb\lambda$,
which is an oscillatory function of order $\mu$,
is bounded by $O(\lambda\sp{\mu+\kappa'})$,
uniformly in $x$ and $y$.
It follows that for
the action
\[
\euf{F}:\; L\sp 1(Y) \to L\sp\infty (X)
\]
to be continuous
we
need to have
$\euf{F}\in I\sp\mu$ with $\mu<-\frac{n+1}{2}-\kappa$.
Thus, in the presence of caustics,
the $L\sp 1\to L\sp\infty$ estimates deteriorate.
On the other hand,
if we assume that $\eub C$ is a local graph,
the mappings
\begin{eqnarray}
&&\euf{F}:\;L\sp 2(Y)\to L\sp 2(X)
\qquad{\rm if}
\quad \mu\le 0,
\label{l2-l2}
\\
&&\euf{F}:\;\hardy\sp 1(Y)\to L\sp 1(X)
\qquad{\rm if}
\quad \mu\le-\frac{n-1}2,
\label{l1-l1}
\\
&&\euf{F}:\;\hardy\sp 1(Y)\to L\sp 2(X)
\qquad{\rm if}
\quad \mu\le-\frac n 2,
\label{l1-l2}
\end{eqnarray}
are continuous,
independently of the presence of caustics.
(The estimate (\ref{l2-l2}) is the classical
$L\sp 2$ bound on Fourier integral operators,
(\ref{l1-l1}) is proved in \cite{MR92g:35252},
and (\ref{l1-l2}) follows from
the $\hardy\sp 1\to L\sp\infty$ estimate on
$\euf{F}\euf{F}\sp\ast$,
which is a pseudodifferential operator
of order $2\mu$.)
The $L\sp p\to L\sp q$ estimates
for $1<p\le q\le 2$
(obtained by interpolation of
(\ref{l2-l2}), (\ref{l1-l1}), and (\ref{l1-l2}))
are also caustic-insensitive.
By duality considerations, the same is true
for $2\le p\le q<\infty$.
We are going to show that for
$1<p\le 2\le q<\infty$
away from a certain neighborhood of the line $q=p'$ the
$L\sp p\to L\sp q$ estimates
are also caustic-insensitive.
In this paper, we only consider the situation when
$\rank d\pi\sb{X\times Y}\ge 2n-2$.

\bigskip

Let us give a short account of our methods.
Let $\eub C\subset T\sp\ast(X)\backslash 0\times T\sp\ast(Y)\backslash 0$
be a canonical relation
which is a local graph.
Let
$\euf{F}=\sum\sb{\lambda}\euf{F}\sb\lambda+\euf{F}\sb{0}$,
$\lambda=2^k$, $k\in\N$
be a Littlewood-Paley decomposition of $\euf{F}\in I\sp\mu(X,Y,\eub C)$.
If the caustic set $\Sigma(\eub C)$ is empty
($\rank d\pi\sb{X\times Y}=2n-1$ everywhere),
then, representing $\euf{F}$ with $\theta\in\R\sp 1$,
one easily checks that
\begin{equation}
\norm{\lambda\sp{\mu+\frac{n+1}{2}}\euf{F}\sb\lambda}
\sb{L\sp 1\to L\sp\infty}<C,
\end{equation}
uniformly in $\lambda$.
Now let $\Sigma(\eub C)\neq \emptyset$,
and assume that
$\rank d\pi\sb{X\times Y}=2n-2$
at $\Sigma(\eub C)$.
Choosing local coordinates $\alpha$
on the unit sphere in $\R^N$,
we introduce a function
\begin{equation}\label{distance-to-caustic-set}
\mathcal{D}=\det\sb{ij}(\abs{\theta}\sp{-1}\phi\sb{\alpha\sb i\alpha\sb j}''),
\qquad
1\le i,\,j\le N-1,
\end{equation}
which measures the
distance to the caustic set.
This function is defined up to a nonzero factor,
which depends on the local coordinates.
We decompose $\euf{F}\sb\lambda$ into
$\sum\sb\sigma\euf{F}\sb{\lambda,\sigma}+\euf{F}\sb{\lambda,\rm nice}$,
where $\sigma=2\sp{-j}$, $j\in\N$,
with respect to the values of $\mathcal{D}$,
so that the Schwartz kernel of $\euf{F}\sb{\lambda,\sigma}$
is localized to the set $\sigma/2\le\abs{\mathcal{D}}\le 2\sigma$
near $\Sigma(\eub C)$,
while the Schwartz kernel of $\euf{F}\sb{\lambda,\rm nice}$
is localized away from $\Sigma(\eub C)$.
When approaching the caustic set,
the $L\sp 1\to L\sp\infty$ estimates
become worse:
\begin{equation}\label{sigma-l1-inf}
\norm{\lambda\sp{-\mu-\frac{n+1}{2}}\euf{F}\sb{\lambda,\sigma}}
\sb{L\sp 1\to L\sp\infty}
\sim\sigma\sp{-1/2}.
\end{equation}
This is the optimal estimate which one expects
from the application of the stationary phase method.
On the other hand,
the $L\sp 1\to L\sp 2$ action of
$\euf{F}\sb{\lambda,\sigma}$ improves near $\Sigma(\eub C)$:
\begin{equation}\label{sigma-l1-l2}
\norm{\lambda\sp{-\mu-\frac{n}{2}}\euf{F}\sb{\lambda,\sigma}}
\sb{L\sp 1\to L\sp 2}
\sim\sigma\sp{\frac 1{2m}}.
\end{equation}
Similarly to the idea from \cite{MR81d:42029},
this estimate is essentially the
``square root'' of the
estimate on the $L\sp 1\to L\sp\infty$ action
of $\lambda\sp{-2\mu-n}\euf{F}\sb{\lambda,\sigma}\sp\ast\euf{F}\sb{\lambda,\sigma}$.
The Schwartz kernel of
$\lambda\sp{-2\mu-n}
\euf{F}\sb{\lambda,\sigma}\sp\ast\euf{F}\sb{\lambda,\sigma}$
is bounded uniformly in $x$, $y$, and $\lambda$,
and hence this operator is bounded from $L\sp 1$ to $L\sp\infty$
(uniformly in $\lambda$).
Moreover,
the Schwartz kernel
involves an inert integration in $\theta$,
and if $\mathcal{D}$ vanishes of order $m$
with respect to $\theta$ 
(as for the caustics of the type $A\sb{m+1}$),
then there is an improvement
$\norm{\lambda\sp{-2\mu-n}\euf{F}\sb{\lambda,\sigma}\sp\ast\euf{F}\sb{\lambda,\sigma}}
\sb{L\sp 1\to L\sp\infty}\sim\sigma\sp{1/m}$
for small values of $\sigma$,
which leads to (\ref{sigma-l1-l2}).

The estimates
(\ref{sigma-l1-inf}) and (\ref{sigma-l1-l2})
allow us to prove that there is some $q\sb c>2$
such that for $2\le q< q\sb c$
the $L\sp 1\to L\sp q$
regularity of $\euf{F}\sb\lambda$ is not affected by caustics.
In essence,
this situation is expressed
by the following obvious lemma:
\begin{lemma}
Let
$B\sb{t}$, $0\le t\le 1$, be a family of complete Banach
spaces
and that $B\sb t\subset B\sb{t'}$,
$\norm{\cdot}\sb{B\sb t}\ge\norm{\cdot}\sb{B\sb{t'}}$
if $t\ge t'$.
Moreover, assume that for any $\xi\in B\sb 1\subset B\sb 0$,
\[
\norm{\xi}\sb{B\sb t}\le\norm{\xi}\sb{B\sb 0}\sp{1-t}\norm{\xi}\sb{B\sb 1}\sp{t}.
\]
Let $\xi\sb j\in B\sb 1$, $j\in\N$,
be a sequence such that
$\norm{\xi\sb j}\sb{B\sb 0}\le a^j$, $0<a<1$,
and that
$\norm{\xi\sb j}\sb{B\sb 1}\le b^j$, $b>1$.
Then
$\sum\sb{j\in\N}\xi\sb j$
converges in $B\sb t$
for $0\le t<t\sb c\equiv\frac{-\ln a}{-\ln a+\ln b}$.
\end{lemma}

This already allows us to calculate the critical value $q\sb c$.
The estimate
(\ref{sigma-l1-inf})
(mapping to $L\sp \infty$) blows up as
$\sigma\sp{-1/2}$,
while the estimate (\ref{sigma-l1-l2})
(mapping to $L\sp 2$) improves as
$\sigma\sp{1/(2m)}$.
Interpolation shows that the bound on
the $L\sp 1\to L\sp q$ mapping
of $\lambda\sp{-\mu-(n+1)/2+1/q}\euf{F}\sb{\lambda,\sigma}$
behaves as $\sigma\sp{-1/2+(1+1/m)/q}$,
which improves for small $\sigma$
if $q<q\sb c=2+2/m$,
so that the mapping
\[
\lambda\sp{-\mu-(n+1)/2+1/q}\euf{F}\sb\lambda:\;L\sp 1(Y)\to L\sp q(X),
\qquad 2\le q<q\sb c,
\]
is bounded uniformly in $\lambda$
and is not affected by the caustics.
These estimates could be interpolated with
the (caustic-insensitive) $L\sp 2\sb{\mu}\to L\sp 2$ estimates
on $\euf{F}\sb\lambda$.
For $1<p\le 2\le q<\infty$,
Littlewood-Paley theory implies that
the mapping $\euf{F}:\;L\sp p\sb{\alpha}\to L\sp q$
is continuous if the mappings
$\lambda^{\alpha}\euf{F}\sb\lambda:\;L\sp p\to L\sp q$
are bounded uniformly in $\lambda$.
This, together with duality considerations,
yields the range of $p$ and $q$
such that the $L\sp p\to L\sp q$ estimates are caustic-insensitive.

\begin{remark}
{\rm
This situation is similar to the $L\sp p\to L\sp p$
regularity of Fourier integral operators
associated to degenerate canonical relations,
when the projections from the canonical relation
are allowed to have singularities.
For example,
as we know from \cite{MR86m:35095},
if both projections
$\eub C\to T\sp\ast(X)$,
$\eub C\to T\sp\ast(Y)$
have at most Whitney fold singularities,
then the operator $\euf{F}\in I^\mu(X,Y,\eub C)$
loses $1/6$ of a derivative in the Sobolev spaces:
$\euf{F}:H\sp{\alpha}(Y)\to H\sp{\alpha-\mu-\frac 1 6}(X)$,
but according to \cite{MR95c:35048}
the $L\sp p\to L\sp p$ regularity of such an operator
for $p\notin (1,3/2)\cup(3,\infty)$
is the same as for
operators associated to local graphs \cite{MR92g:35252}:
$\euf{F}:L\sp p\sb{\alpha}(Y)
\to L\sp p\sb{\alpha-\mu-(n-1)\abs{\delta\sb{p}}}(X)$,
$\delta\sb{p}=\frac 1 p -\frac 1 2$.
Similar estimates on operators with
one-sided Whitney folds were derived in \cite{hardy}.
This time, one uses
the Phong-Stein decomposition \cite{MR93g:58144}
of $\euf{F}$
with respect to the distance to the critical variety
$\Sigma$
(where the projections from the canonical
relation become singular).
This distance is measured by the function
\[
h=
\abs{\theta}^{N-n}
\det\left[
\begin{array}{cc}
\phi\sb{xy}''&\phi\sb{x\theta}''
\\
\phi\sb{\theta y}''&\phi\sb{\theta\theta}''
\end{array}
\right],
\]
which is proportional to the determinants
of the Jacobi matrices of projections
from $\eub C$.
(The factor in the definition of $h$
is chosen so that $h$ is homogeneous of degree $0$
in $\theta$.)
We decompose
$\euf{F}=\sum\sb\hvar \euf{F}\sb\hvar+\euf{F}\sb{\rm smooth}$,
where $\hvar=2\sp{-j}$, $j\in\N$.
The operator $\euf{F}\sb\hvar$ is
obtained from $\euf{F}$ by
localizing its integral kernel
to the variety where
$\hvar/2\le \abs{h}\le 2\hvar$,
and the projections from $\eub C$
have no singularities on the support
of the integral kernel of $\euf{F}\sb{\rm smooth}$.
The main observation is that
while near $\Sigma$
the Sobolev estimates become worse,
$
\norm{\euf{F}\sb\hvar}\sb{L\sp 2\sb{\mu}\to L\sp 2}\sim \hvar\sp{-1/2},
$
the Hardy space to $L\sp 1$ estimates improve
due to smaller size of the support of the integral kernel:
$
\norm{\euf{F}\sb\hvar}\sb{\hardy\sp 1\sb{\mu+\frac{n-1}{2}}\to L\sp 1}
\sim \hvar.
$
}
\end{remark}

\bigskip

Caustics of Lagrangian distributions
are discussed in Section \ref{sect-caustics}.
The main results
(Theorems \ref{theorem-index-m-lp-lq} and \ref{theorem-index-m-hardy})
are stated in Section \ref{sect-main-results}.
The proof of $L\sp p\to L\sp q$ estimates
is in Section \ref{sect-microlocal}.
The sharp $\hardy\sp 1\to L\sp q$ estimates are proved
in Section \ref{sect-microlocal-hardy}.
We apply our results
to the estimates on the half-wave operator
in Section \ref{sect-geodesic-flow}.

The consistency of definition 
(\ref{distance-to-caustic-set})
of the distance $\mathcal{D}$ to the caustic set
is proved in Appendix \ref{appendix-consistency}.
The technical lemma
($\hardy\sp 1\to L\sp\infty$ bounds on pieces)
which allows us to obtain
$\hardy\sp 1\to L\sp q$ estimates
is proved in Appendix \ref{appendix-hardy}.

\sect{Caustics of Lagrangian distributions}\label{sect-caustics}

\subsection{Symbols}

We will use the class of classical (polyhomogeneous) symbols,
in the sense of \cite{MR95h:35255}.

\begin{definition}\label{def-classical-symbol}
A smooth function $a(x,\theta)$ on $X\times \R^N$
is called a symbol of order $d$ in $\theta$
if for any multi-indices $\alpha\in\Z\sb{+}^n$
and $\beta\in\Z\sb{+}^N$
\[
\abs{\p\sb{x}^\alpha\p\sb\theta^\beta a(x,\theta)}
\le C\sb{\alpha,\beta}(1+\abs{\theta})\sp{d-\abs{\beta}},
\qquad
{\rm for\ all}
\quad
(x,\theta)\in X\times\R^N,
\]
where
$\abs{\beta}=\beta\sb 1+\dots+\beta\sb N$.

We denote the class of symbols of order $d$ by $S^d(X\times\R^N)$
or simply by $S^d$.

The class of classical (or polyhomogeneous) symbols
$S\sb{\rm cl}^d(X\times\R^N)$
consists of symbols $a(x,\theta)\in S^d(X\times\R^N)$
that satisfy an asymptotic development of the form
\begin{equation}\label{asymptotic-development}
a(x,\theta)
\sim\sum\sb{j=0}\sp\infty a\sb j(x,\theta),
\end{equation}
where $a\sb j$ are smooth functions on $X\times\R^N$
positively homogeneous of degree $d-j$ for $\abs\theta\ge 1$:
\begin{equation}
a\sb j(x,\tau\theta)=\tau\sp{d-j}a\sb j(x,\theta)
\quad{\rm if}\quad \abs\theta\ge 1, \quad \tau\ge 1.
\end{equation}
\end{definition}

The asymptotic development (\ref{asymptotic-development})
means that we have
\begin{equation}
a(x,\theta)
-\sum\sb{j=0}\sp{k-1}a\sb j(x,\theta)=O(\abs\theta\sp{d-k})
\quad{\rm for}\quad \abs\theta\ge 1
\end{equation}
and similar estimates for the derivatives.

\subsection{Oscillatory functions}

Let $X$ be a $C\sp\infty$ manifold
and $\Lambda\subset T\sp\ast(X)$
a $C\sp\infty$ Lagrangian submanifold.
We say that $\psi\in C\sp\infty(X\times\R^k)$
parametrizes $\Lambda$ (locally) if
\begin{equation}\label{psi-non-degenerate}
d\sb{(x,\alpha)}d\sb\alpha\psi
\quad{\rm has\ rank\ }k{\rm\ when\ }
d\sb\alpha\psi=0
\end{equation}
and $\Lambda$ is locally given by
\begin{equation}\label{lambda-psi}
\Lambda\sb\psi=\{
(x,d\sb{x}\psi(x,\alpha))\sothat d\sb\alpha\psi(x,\alpha)=0
\}.
\end{equation}

\begin{definition}\label{def-osc-function}
Let $\Lambda$ be a $C\sp\infty$ Lagrange submanifold
in $T\sp\ast(X)$.
An oscillatory function $u(x,\tau)$ of order $\mu$
defined by $\Lambda$
is a locally finite (in $X$) sum of integrals of the form
\[
I(x,\tau)=\int e\sp{i\tau\psi(x,\tau,\alpha)}b(x,\tau,\alpha)\,d\alpha,
\]
where $\alpha=(\alpha\sb 1,\dots,\alpha\sb k)$,
$k\in\N$,
$\psi$ satisfies (\ref{psi-non-degenerate}),
$\Lambda\sb\psi$ is a piece of $\Lambda$
and $b(x,\tau,\alpha)\in S\sb{\rm cl}\sp{\mu+\frac k 2}$
(a classical symbol of order $\mu+\frac k 2$ in $\tau$)
vanishes for $\alpha$ outside a fixed compact set in $\R^k$.
\end{definition}

\begin{definition}[Duistermaat \cite{MR96m:58245}]\label{def-caustics}
Let $i:\Lambda\to T\sp\ast(X)$
be an immersed Lagrange manifold in $T\sp\ast(X)$.
The caustic $c(\Lambda)$ of $\Lambda$
is the projection into $X$
of the set $\Sigma(\Lambda)$ of points in $i(\Lambda)$
where $i$ is not transversal to the fibers.
At each point $x\sb 0\in X$ the order of the caustic
is defined as the infimum $\kappa(x\sb 0)$
of the numbers $\kappa'$
such that $u(x,\tau)=O(\tau\sp{\mu+\kappa'})$
for $\tau\to\infty$,
uniformly for $x$ in a neighborhood of $x\sb 0$,
for any oscillatory function $u$ of order $\mu$
defined by $\Lambda$.
\end{definition}
Of course $\kappa(x\sb 0)=0$
for $x\sb 0\in\pi(\Lambda)\backslash c(\Lambda)$
and $\kappa(x\sb 0)\le k/2$
where $k$ is the maximum of the dimensions
of the intersections
$T\sb{(x\sb 0,\xi)}(\Lambda)\cap T\sb{(x\sb 0,\xi)}({\rm fiber})$
where $(x\sb 0,\xi)\in\Lambda$.

\subsection{Lagrangian distributions}

We also need to define caustics of the conic Lagrangian submanifolds.
Let $X$ be a $C\sp\infty$ manifold
and $\Lambda\subset T\sp\ast(X)$
be a conic $C\sp\infty$ Lagrangian submanifold.
We say that $\phi\in C\sp\infty(X\times\R^N)$
parametrizes $\Lambda$ (locally) if
\begin{equation}\label{phi-non-degenerate}
d\sb{(x,\theta)}d\sb\theta\phi
\quad{\rm has\ rank\ }N{\rm\ when\ }
d\sb\theta\phi=0
\end{equation}
and $\Lambda$ is locally given by
\begin{equation}\label{lambda-phi}
\Lambda\sb\phi=\{
(x,d\sb{x}\phi(x,\theta))\sothat d\sb\theta\phi(x,\theta)=0
\}.
\end{equation}

\begin{definition}[H\"ormander]\label{def-lagrangian-density}
Let $\Lambda$ be a $C\sp\infty$ conic Lagrangian submanifold
in $T\sp\ast(X)$.
A Lagrangian distribution $u(x)$ of order $\mu$
defined by $\Lambda$
is a locally finite (in $X$) sum of integrals of the form
\[
u(x)=\int e\sp{i\phi(x,\theta)}a(x,\theta)\,d\theta,
\]
where $\theta\in\R^N$,
$\phi$ satisfies (\ref{phi-non-degenerate}),
$\Lambda\sb\phi$ is a piece of $\Lambda$
and 
\[
a(x,\theta)\in S\sp{\mu-\frac N 2+\frac {\dim X} 4}(X\times\R^N).
\]
\end{definition}

We pick a smooth function
$\rho\in C\sp\infty\sb{0}([-2,2])$,
$\rho\ge 0$, $\rho\at{[-1,1]}\equiv 1$.
Define $\beta(t)=\rho(t)-\rho(2t)$
for $t>0$, $\beta\equiv 0$ for $t\le 0$,
so that $\beta\in C\sp\infty\sb{0}([\frac 1 2,2])$.
We introduce the Littlewood-Paley decomposition
of $u(x)$:
\begin{equation}
u\sb\lambda(x)=\int e\sp{i\phi(x,\theta)}\beta(\abs\theta/\lambda)
a(x,\theta)\,d\theta.
\end{equation}

\begin{definition}\label{def-conic-caustics}
Let $\Lambda$
be a conic Lagrangian manifold in $T\sp\ast(X)$.
The caustic $c(\Lambda)$ of $\Lambda$
is the projection into $X$
of the set $\Sigma(\Lambda)$ of points in $\Lambda$
where $\rank d\pi\sb{X}<\dim X-1$.
At each point $x\sb 0\in X$ the order of the caustic
is defined as the infimum $\kappa(x\sb 0)$
of the numbers $\kappa'$
such that 
$u\sb\lambda(x)=O(\lambda\sp{\mu+\frac{\dim X} 4 +\frac 1 2 +\kappa'})$
for $\lambda\to\infty$,
uniformly for $x$ in a neighborhood of $x\sb 0$,
for any Lagrangian distribution $u$ of order $\mu$
defined by $\Lambda$.
\end{definition}

\begin{definition}
We say that $\Lambda$ has a caustic of corank $k\ge 1$
at a point $p\sb 0\in\Sigma(\Lambda)$
if
\[
\rank d\pi\sb{X}\at{p\sb 0}=\dim X-1-k.
\]
\end{definition}

\begin{lemma}\label{lemma-consistency}
Let $\Lambda$ be a smooth closed conic Lagrangian submanifold
of $T\sp\ast(X)\backslash 0$.
Let $\phi(x,\theta)\in C\sp\infty(X\times\R^N)$
be a smooth non-degenerate phase function
which parametrizes $\Lambda$:
\[
\Lambda=\{(x,d\sb{x}\phi(x,\theta))
\sothat d\sb\theta\phi(x,\theta)=0\}.
\]
Let $\alpha=\{\alpha\sb i\}$, $1\le i\le N-1$,
be local coordinates on the unit sphere $\S\sp{N-1}$.
We use $(\lambda,\alpha)\in \R\sb{+}\times \S\sp{N-1}$
as local coordinates in $\R^N$.
Then
$\mathcal{D}
=\det\sb{ij}(\lambda\sp{-1}\phi\sb{\alpha\sb i\alpha\sb j}''\at{\Lambda})$,
$1\le i,\,j\le N-1$,
is a smooth function on $\Lambda$ defined up to a nonzero factor:
\[
\mathcal{D}\in
C\sp\infty(\Lambda)\slash C\sp\infty\sb\times(\Lambda).
\]
\end{lemma}
This statement is intuitively clear, since
$\mathcal{D}$ vanishes precisely on the caustic set
$\Sigma(\Lambda)$
where the rank of the projection from $\Lambda$ onto
$X$ is smaller than $\dim X-1$.
Still,
since we need to know that the order of vanishing
of $\mathcal{D}$ at $\Sigma(\Lambda)$
in particular directions
does not depend on the number of oscillatory variables
and the choice of local coordinates,
we will give a detailed argument in Appendix \ref{appendix-consistency}.

\begin{definition}\label{def-simple-caustics}
We say that the caustic at a point $p\sb 0\in\Sigma(\Lambda)$ is simple 
if it is of corank $k=1$,
so that 
\[
\rank d\pi\sb{X}\at{p\sb 0}=\dim X-1-k=\dim X-2,
\]
and if $d\sb{(x,\theta)}\mathcal{D}\at{p\sb 0}\ne 0$.
\end{definition}

\begin{definition}\label{def-index-of-caustic}
We say that the simple caustic at a point $p\sb 0\in\Sigma(\Lambda)$ 
is of index $m\in\N$
if $m$ is the smallest integer so that
there exists a vector field
$V\in C\sp\infty(\Gamma(T(\Lambda)))$,
$V\at{\Sigma(\Lambda)}\in\ker d\pi\sb{X}$,
such that
\[
V\sp{m}\mathcal{D}(p\sb 0)\ne 0,
\]
where
\[
\mathcal{D}=\det\sb{ij}(\abs{\theta}\sp{-1}\phi\sb{\alpha\sb i\alpha\sb j}''),
\qquad 1\le i,\,j\le N-1.
\]
\end{definition}

\begin{remark}
Consistency of this definition
(independence of the choice of the phase function $\phi$
that parametrizes $\Lambda$
and independence of the choice of local coordinates $\alpha$ on $\mathbb{S}\sp{N-1}$)
follows from Lemma \ref{lemma-consistency}.
\end{remark}

\begin{example}
Let $\theta\in\R\sp 2$ and $\lambda=\abs{\theta}$.
Then $\theta/\abs{\theta}\in\mathbb{S}\sp 1$.
Denote a local coordinate on $\mathbb{S}\sp 1$ by $\alpha$.
Consider the phase function
$\phi(x,\theta)=\abs{\theta}\Phi(x,\alpha)$,
with
\[
\Phi(x,\alpha)=\alpha\sp{m+2}+x\sb 1\alpha\sp m+\dots + x\sb{m}\alpha+x\sb{m+1}.
\]
This is the model example of a caustic of the type $A\sb{m+1}$.
The corresponding Lagrangian is
$
\Lambda=\{x,d\phi(x,\theta)\sothat\phi\sb\theta'=0\},
$
which can be written as
\[
\Lambda=\{x,\lambda d\sb x\Phi(x,\alpha)
\sothat \Phi(x,\alpha)=0,\ \Phi\sb\alpha'(x,\alpha)=0\}.
\]
The Lagrangian could be parametrized by $(x',\lambda,\alpha)$,
where $x'=(x\sb 1,\dots,x\sb{m-1})$.
At $\Lambda$, 
one can express $x\sb m$ and $x\sb{m+1}$ as
$x\sb m=X\sb m(x',\alpha)$,
$x\sb{m+1}=X\sb{m+1}(x',\alpha)$.

The kernel of the map $\pi:\;\Lambda\to X$ always contains 
the tangent vector $\p/\p\lambda$.
On the caustic set
\[
c(\Lambda)=\{(x',\alpha)\in\Lambda\sothat \mathcal{D}(x',\alpha)
=\Phi\sb{\alpha\alpha}''(x',x\sb m(x',\alpha),x\sb{m+1}(x',\alpha),\alpha)=0\}
\]
one also has $\p/\p\alpha\in\ker d\pi$.
Since $1\le \dim\ker d\pi\le 2$
and $d\mathcal{D}=d(\Phi\sb{\alpha\alpha}''\at{\Lambda})\ne 0$,
the caustics are simple
in the sense of Definition \ref{def-simple-caustics}.

Consider the vector field $V=\p/\p\alpha\in C\sp\infty(\Gamma(T\Lambda))$,
$V\at{c(\Lambda)}\in\ker d\pi$.
Since
\[
V\sp{m}\mathcal{D}(x',\alpha)
=\p\sb\alpha \sp{m}(\Phi\sb{\alpha\alpha}''\at{\Lambda})
=\p\sb\alpha \sp{m}\Phi\sb{\alpha\alpha}''(x',x\sb m(x',\alpha),x\sb{m+1}(x',\alpha),\alpha)
\ne 0,
\]
one concludes that the caustic is of index (at most) $m$.

\end{example}

\begin{remark}
While caustics of the type $A\sb{m+1}$, $m\ge 1$, correspond to simple caustics of index $m$,
the converse is not necessarily true,
except when $m=1$ and $2$.
\end{remark}

Let us show how to prove that
simple caustics of index $m=1$ and $2$
necessarily correspond to caustics of the type $A\sb{2}$ and $A\sb{3}$,
respectively.
We first reduce the number of oscillatory variables to $N=2$,
denote by $\alpha$ the local coordinate on $\mathbb{S}\sp 1$,
and define
\[
\Phi(x,\alpha)=\phi(x,\theta/\abs{\theta})=\phi(x,\theta)/\abs{\theta}.
\]
It suffices to notice that
$\phi$ has 
the caustic of the type $A\sb{m+1}$
at the point $(x\sb 0,\alpha\sb 0)$
if
$\p\sb{\alpha}\sp{j}\Phi(x\sb 0,\alpha\sb 0)=0$,
$j<m+2$, $\p\sb{\alpha}\sp{m+2}\Phi(x,\alpha)\ne 0$,
and if the differentials
\[
d\Phi(x,\alpha),
\quad
d\Phi\sb\alpha'(x,\alpha),
\quad
\dots,
\quad
d\Phi\sb{\alpha\dots}\sp{(m)}(x,\alpha)
\]
are linearly independent.
In the case $m=1$,
the linear independence of 
$d\Phi$ and $d\Phi\sb{\alpha}'$
follows from the non-degeneracy assumption on $\phi$
(the differentials $d\phi\sb{\theta\sb j}$ are linearly independent).

To settle the case $m=2$,
we additionally need to check that 
$d\Phi\sb{\alpha\alpha}''$
is linearly independent of 
$d\Phi$ and $d\Phi\sb\alpha'$.
We only need to notice that two latter differentials
vanish identically on $T\Lambda$,
while the differential
$d\Phi\sb{\alpha\alpha}''=d\mathcal{D}$
was assumed to be different from zero
(see Definition \ref{def-simple-caustics}).

\sect{Main results}\label{sect-main-results}

Consider a Fourier integral operator
\begin{equation}\label{fio}
\euf{F} u(x)=\int\sb{\R^N\times Y}
e\sp{i\phi(x,\theta,y)}a(x,\theta,y)u(y)\,d\theta\,dy,
\end{equation}
where $X$ and $Y$ are two smooth manifolds
(unless stated otherwise, we assume that $\dim X=\dim Y=n$),
$\ a(x,\theta,y)\in S^d\sb{\rm cl}(X\times\R^N\times Y)$
is a symbol of order $d$ in $\theta$,
compactly supported in $X\times Y$.
(We restrict the consideration
to classical (polyhomogeneous) symbols, denoted by $S\sb{\rm cl}$.)
The function $\phi$ is a non-degenerate phase:
the differentials
$d\sb{(x,\theta,y)}\phi\sb{\theta\sb j}'$, $1\le j\le N$,
are linearly independent,
so that
\[
\varSigma\sb\phi=
\{
(x,\theta,y)\sothat \phi\sb\theta'(x,\theta,y)=0
\}
\]
is a smooth submanifold of $X\times \R^N\times Y$
of dimension $\dim X+\dim Y$.
We write $\euf{F}\in I\sp\mu(X,Y,\eub C)$,
where the order of the operator
is defined by $\mu=d+\frac{N}{2}-\frac{\dim X+\dim Y}{4}$
and $\eub C$ is the associated canonical relation:
\[
\varSigma\sb\theta
\stackrel{\cong}\longrightarrow\eub C
\equiv\{(x,d\sb{x}\phi(x,\theta,y)),(y,-d\sb{y}\phi(x,\theta,y))
\sothat \phi\sb\theta'(x,\theta,y)=0\}.
\]
According to \cite{MR52:9299},
the minimal number of oscillatory variables
is equal to 
\[
N\sb{\rm min}=\dim X+\dim Y-\min\rank d\pi\sb{X\times Y}.
\]
If $\eub C$ has non-empty caustic set $\Sigma(\eub C)$,
then $N\sb{\rm min}>1$.

\begin{definition}\label{def-c-caustics}
We will say that the canonical relation
$\eub C\subset T\sp\ast(X)\backslash 0\times T\sp\ast(Y)\backslash 0$
has a caustic at a point 
$p\sb 0=((x\sb 0,\xi\sb 0),(y\sb 0,\eta\sb 0))\in \eub C$
if the twisted canonical relation
$\eub C'=\{((x,\xi),(y,-\eta))\sothat((x,\xi),(y,\eta))\in\eub C\}$,
which is a conic Lagrangian submanifold 
of $T\sp\ast(X\times Y)\backslash 0$,
has a caustic at a point 
$((x\sb 0,\xi\sb 0),(y\sb 0,-\eta\sb 0))\subset \eub C' $.
We will not distinguish 
the caustics of $\eub C$ and $\eub C'$.

\end{definition}

The following result
is an immediate consequence of Definitions \ref{def-conic-caustics}
and \ref{def-c-caustics}:
 
\begin{theorem}\label{easy-theorem}
Let $X$, $Y$ be two smooth manifolds (possibly of different dimension),
and let $\eub C\subset T\sp\ast(X)\backslash 0\times T\sp\ast(Y)\backslash 0$
be a smooth canonical relation.
Let the Fourier integral operator
$\euf{F}\in I\sp\mu(X,Y,\eub C)$
have its symbol compactly supported 
in $X\times Y$.
(The symbol of $\euf{F}$ does not have to be polyhomogeneous.)
If
$\eub C$
has caustics of order at most $\kappa$,
then
\begin{equation}
\euf{F}:\;L\sp 1(Y)\to L\sp\infty(X)
\qquad {\rm if} \quad \mu<-\frac{\dim X+\dim Y}{4}-\frac 1 2-\kappa.
\end{equation}

Further, assume that $\dim X=\dim Y=n$
and that $\eub C$ is a local graph.
Then 
\begin{equation}
\euf{F}:\;L\sp p\sb{\mu+(n+1+2\kappa)\delta\sb{p}}(Y)\to L\sp {p'}(X),
\qquad
1<p\le 2,
\quad
\delta\sb p =\frac 1 p- \frac 1 2.
\end{equation}
\end{theorem}

\begin{proof}
The first part of the theorem follows from the trivial estimate
\begin{equation}\label{trivial-estimate}
\norm{
\euf{F}\sb\lambda u}
\sb{L\sp\infty}\le C
\lambda^{\mu+\frac{\dim X+\dim Y}{4}+\frac 1 2+\kappa}
\norm{u}\sb{L\sp 1}.
\end{equation}
For the second part, we interpolate 
(\ref{trivial-estimate})
with
$\norm{\euf{F}\sb\lambda}\sb{L\sp 2\to L\sp 2}\le C\lambda^{\mu}$
and apply Littlewood-Paley theory.
\end{proof}

\begin{definition}
For our convenience, we introduce the map
\begin{equation}
(\,\cdot\,,\cdot\,)\sp\dagger:\;(p,q)\mapsto(p,q)\sp\dagger=(1/p,1/q).
\end{equation}

\end{definition}

\begin{figure}[htbp]
\input lplq-caustics-regions.tex
\caption{Regions
$\euf{A}\sb m$
$\euf{B}\sb m$,
and
$\euf{C}\sb m$
in $(1/p,1/q)$-plane.
}
\label{fig-omega-regions}
\end{figure}

\begin{definition}\label{def-omega-regions}
We define
\[
p\sb m=2-\frac{2}{m+2},
\qquad
\ q\sb m=p\sb m'=2+\frac{2}{m}
\]
and introduce the following regions in $(1/p,1/q)$-plane
(see Figure \ref{fig-omega-regions}):

$\euf{A}\sb m$
is the open triangle with the vertices at
$(2,2)\sp\dagger$, $(1,1)\sp\dagger$, and $(1,q\sb m)\sp\dagger$.

$\euf{B}\sb m$
is the open triangle with the vertices at
$(\infty,\infty)\sp\dagger$, $(2,2)\sp\dagger$, and $(p\sb m,\infty)\sp\dagger$.

$\euf{C}\sb m$
is the open convex hull of points
$(1,\infty)\sp\dagger$, $(1,q\sb m)\sp\dagger$, $(2,2)\sp\dagger$,
and $(p\sb m,\infty)\sp\dagger$.
\end{definition}

The following is our main result:

\begin{theorem}\label{theorem-index-m-lp-lq}
Let $X$, $Y$ be two smooth manifolds,
$\dim X=\dim Y=n$,
and let $\eub C\subset T\sp\ast(X)\backslash 0\times T\sp\ast(Y)\backslash 0$
be a smooth canonical relation which is a local graph.
Assume that
$\eub C$ has only simple caustics of index at most $m$, $m\in\N$.
Let $\euf{F}\in I\sp\mu(X,Y,\eub C)$
have the classical (polyhomogeneous) symbol
with compact support in $X\times Y$.
Then for
$(p,q)\sp\dagger\in\euf{A}\sb m
\cup\euf{B}\sb m$
the $L\sp p\to L\sp q$ estimates on $\euf F$ are caustics-insensitive.

Precisely,
\begin{eqnarray}
&&\euf{F}
:\;
L\sp p\sb{\mu+n\delta\sb{p}+\delta\sb{q}}(Y)
\to L\sp q(X),
\qquad (p,q)\sp\dagger\in\euf{A}\sb m,
\\
&&\euf{F}
:\;
L\sp p\sb{\mu+n\delta\sb{q}+\delta\sb{p}}(Y)
\to L\sp q(X),
\qquad (p,q)\sp\dagger\in\euf{B}\sb m,
\end{eqnarray}
where
\[
\delta\sb{p}=\frac 1 p-\frac 1 2,
\qquad
\delta\sb{q}=\frac 1 2 -\frac 1 q.
\]

\noindent
For $(p,q)\sp\dagger\in\euf{C}\sb m$,
the estimates depend on the order of the caustic,
given by $\kappa=\frac 1 2 -\frac{1}{m+2}$:
\begin{eqnarray}
&&\euf{F}
:\;
L\sp p\sb{\mu
+n\delta\sb{p}
+(\delta\sb{p}+\delta\sb{q})(1/2+\kappa)
+(\delta\sb{q}-\delta\sb{p})
}(Y)\to L\sp q(X),
\qquad
(p,q)\sp\dagger\in\euf{C}\sb m,
\quad q\le p',
\\
&&\euf{F}:\;
L\sp p\sb{\mu
+n\delta\sb{q}
+(\delta\sb{p}+\delta\sb{q})(1/2+\kappa)
+(\delta\sb{p}-\delta\sb{q})
}(Y)\to L\sp q(X),
\qquad
(p,q)\sp\dagger\in\euf{C}\sb m,
\quad q>p'.
\end{eqnarray}
\end{theorem}

\begin{remark}
We restrict the consideration to the class
of classical symbols
$S\sp{d}\sb{\rm cl}\subset S\sp{d}\sb{1,0}$
in order to simplify
the proof of Lemma \ref{lemma-i-bound}.
\end{remark}

\begin{remark}
Note that $\kappa=\frac 1 2 -\frac{1}{m+2}$
is the order of a caustic of the type $A\sb{m+1}$.
\end{remark}

\begin{remark}
We need the assumption that
$\eub C$ is a local graph
in order to interpolate with the $L\sp 2$-based
Sobolev estimates on $\euf{F}$:
\begin{equation}\label{l2-continuity}
\euf{F}:\;L\sp 2\sb\mu(Y)\to L\sp 2(X).
\end{equation}
The argument could immediately be adapted to the case
when the projections from $\eub C$ have singularities,
as long as $\eub C\to X$ and $\ \eub C\to Y$ are assumed to be submersions.
In this case, one only needs to modify (\ref{l2-continuity}),
taking into account the loss of derivatives
due to singularities of the projections; see \cite{MR1964817}.

\end{remark}

\begin{remark}
The sharp estimates on the line $p=q$
follow from \cite{MR92g:35252}:
\[
\euf{F}:\;\hardy\sp 1\sb{\mu+\frac{n-1}{2}}(Y)\to L\sp 1 (X).
\]
This map can be interpolated with the
continuous $L\sp 2\sb\mu\to L\sp 2$ action.
(Generalization for operators with degenerate canonical relations
is obtained in \cite{hardy}.)
\end{remark}

\begin{remark}
On the line segments $((2,2)\sp\dagger,(1,q\sb m)\sp\dagger)$,
$((2,2)\sp\dagger,(p\sb m,\infty)\sp\dagger)$,
and on the lines $p=1$ and $q=\infty$
the stated estimates hold with the loss of $\epsilon>0$.
\end{remark}

In particular cases, we also have
sharp $\hardy\sp 1\to L\sp q$ and
$L\sp p\to \bmo$ estimates,
as stated in the next theorem.

\begin{theorem}\label{theorem-index-m-hardy}
Let $X$, $Y$ be two smooth manifolds,
$\dim X=\dim Y=n$,
and let $\eub C\subset T\sp\ast(X)\backslash 0\times T\sp\ast(Y)\backslash 0$
be a smooth canonical relation
such that both $\eub C\to X$ and $\eub C\to Y$ are submersions.
Assume that $\eub C$ has only caustics of the type $A\sb{m+1}$ with $m=1$ or $2$.
Let $\euf{F}\in I\sp\mu(X,Y,\eub C)$
have the classical (polyhomogeneous) symbol
compactly supported in $X\times Y$.
Then
\begin{eqnarray}
&&\euf{F}:\;\hardy\sp 1\sb{\mu+\frac n 2 +\delta\sb{q}}(Y)\to L\sp q
(X),
\qquad 2\le q<q\sb m,
\\
&&\euf{F}:\;\hardy\sp 1
\sb{\mu+\frac n 2 +\delta\sb{q}
+\kappa\frac{\delta\sb{q}-\delta\sb{q\sb m}}{1/2-\delta\sb{q\sb m}}
}(Y)\to L\sp q(X),
\qquad q\sb m<q\le\infty.
\end{eqnarray}

\end{theorem}

\begin{remark}
In this theorem,
we do not need $\eub C$ to be a local graph.
\end{remark}

\begin{remark}
$L\sp p(Y)\to \bmo(X)$ estimates on $\euf{F}$
for $1<p<2$
are obtained by duality.
Other $L\sp p\to L\sp q$ estimates
can be obtained by interpolation
with $L\sp 2$ Sobolev estimates.
\end{remark}

\sect{
Microlocal techniques: decompositions and interpolations
}\label{sect-microlocal}

In this section, we prove Theorem \ref{theorem-index-m-lp-lq}.

\subsection{Dyadic decompositions}
We pick a smooth function
$\rho\in C\sp\infty\sb{0}([-2,2])$,
$\rho\ge 0$, $\rho\at{[-1,1]}\equiv 1$.
Define $\beta\in C\sp\infty\sb{0}([\frac 1 2,2])$
by $\beta(t)=\rho(t)-\rho(2t)$
for $t>0$, $\beta\equiv 0$ for $t\le 0$.
The functions
$\rho$ and $\beta$ define dyadic partition of unity:
for any $t\in\R$, 
\[
\sum\sb\pm\sum\sb{j\in\N}\beta(\pm 2\sp{-j}t/2)
+\rho(\abs{t})
=1.
\]
We use the partition of unity
which is the Littlewood-Paley decomposition with respect to the magnitude of $\abs\theta$
and the dyadic decomposition with respect to the distance $\mathcal{D}$ from $\Sigma(\eub C)$:
\[
1=\left(
\sum\sb{\lambda=2^l,\,l\in\N}
\!\!\!
\beta(2^{-l}\abs\theta)
+\rho(\abs\theta)
\right)
\left(
\sum\sb\pm
\sum\sb{j=1}\sp{j\sb 0-1}
\beta(\pm 2\sp{j}\mathcal{D})
+\rho(2\sp{j\sb 0}\abs{\mathcal{D}})
+(1-\rho(2\abs{\mathcal{D}}))
\right).
\]

We define
\begin{equation}\label{f-lambda-sigma}
\euf{F}\sb{\lambda,\pm\sigma} u(x)=
\int\sb{\R^N\times Y}
e\sp{i\phi(x,\theta,y)}
\beta(\pm\mathcal{D}(x,\theta,y)/\sigma)
\beta(\abs{\theta}/\lambda)
a(x,\theta,y)u(y)\,d\theta d y,
\end{equation}
\begin{equation}\label{f-tilde-lambda-sigma}
\tilde{\euf{F}}\sb{\lambda,\sigma} u(x)=
\int\sb{\R^N\times Y}
e\sp{i\phi(x,\theta,y)}
\rho(\mathcal{D}(x,\theta,y)/\sigma)
\beta(\abs{\theta}/\lambda)
a(x,\theta,y)u(y)\,d\theta d y.
\end{equation}
We also define
\begin{equation}
\euf{F}\sb{\rm smooth} u(x)=
\int\sb{\R^N\times Y}
e\sp{i\phi(x,\theta,y)}
\rho(\abs{\theta})
a(x,\theta,y)u(y)\,d\theta d y
\end{equation}
and
\begin{equation}
\euf{F}\sb{\rm nice} u(x)=
\int\sb{\R^N\times Y}
e\sp{i\phi(x,\theta,y)}
\left[
\big(1-\rho(\abs{\theta})\big)\big(1-\rho(2\abs{\mathcal{D}}\big)
\right]
a(x,\theta,y)u(y)\,d\theta d y.
\end{equation}
There is a decomposition
\begin{equation}\label{f-decomposed}
\euf{F}=
\sum\sb\pm
\sum\sb{\lambda}
\sum\sb{\sigma>\sigma\sb 0(\lambda)}
\euf{F}\sb{\lambda,\pm\sigma}
\ +\ 
\sum\sb{\lambda}
\tilde{\euf{F}}\sb{\lambda,\sigma\sb 0(\lambda)}
\ +\ \euf{F}\sb{\rm smooth}
\ +\ \euf{F}\sb{\rm nice},
\end{equation}
where both $\lambda$ and $\sigma$ run over powers of $2$:
\[
\lambda=2\sp l, \quad l\in\N,
\qquad
\sigma=2\sp{-j},\quad 
1\le j<j\sb 0(\lambda)
\equiv[[\log\sb 2\lambda\sp{\frac{m}{m+2}}]].
\]
We set $\sigma\sb 0(\lambda)=2\sp{-j\sb 0(\lambda)}$, 
so that
\begin{equation}\label{sigma-0}
\sigma\sb 0(\lambda)\approx\lambda\sp{-\frac{m}{m+2}}.
\end{equation}
We use the symbol ``$\approx$'' to indicate 
that the quantities differ at most by a factor of $2$.

The operator $\euf{F}\sb{\rm smooth}$ 
is infinitely smoothing and can be discarded.
Since there are no caustics on the support of
$1-\rho(2\abs{\mathcal{D}})$,
the operator $\euf{F}\sb{\rm nice}$ can also be discarded.
The estimates on operators
$\euf{F}\sb{\lambda,\pm\sigma}$
are the same independent of the sign,
and the treatment is the same;
we will only consider the ``$+$''-case.

\subsection{$L\sp 1 \to L\sp\infty$ estimates}

\begin{proposition}\label{prop-l1-inf}
Let $\eub C$ and $\euf{F}\in I\sp\mu(X,Y,\eub C)$
be as in Theorem \ref{theorem-index-m-lp-lq},
and let $\euf{F}\sb{\lambda,\sigma}$,
$\tilde{\euf{F}}\sb{\lambda,\sigma}$
be given by {\rm (\ref{f-lambda-sigma}), (\ref{f-tilde-lambda-sigma}).}
Then
\begin{equation}\label{l1-inf-off-caustics}
\norm{
\euf{F}\sb{\lambda,\sigma}
}\sb{L\sp 1\to L\sp\infty}
\le C\lambda\sp{\mu+\frac{n+1}{2}}\sigma\sp{-\frac 1 2},
\end{equation}
\begin{equation}\label{l1-inf-at-caustics}
\norm{\euf{F}\sb{\lambda,\sigma}}\sb{L\sp 1\to L\sp\infty}
+\norm{\tilde{\euf{F}}\sb{\lambda,\sigma}}\sb{L\sp 1\to L\sp\infty}
\le C\lambda\sp{\mu+\frac{n+2}{2}}\sigma\sp{\frac 1 m}.
\end{equation}
\end{proposition}

\begin{remark}
The value $\sigma\sb 0(\lambda)\approx\lambda\sp{-\frac{m}{m+2}}$
in (\ref{sigma-0})
is chosen so that the estimates
(\ref{l1-inf-off-caustics}),
(\ref{l1-inf-at-caustics})
coincide at $\sigma=\sigma\sb{0}(\lambda)$.
\end{remark}

\begin{proof}
We use the representation of $\euf{F}$
with the minimal possible number of oscillatory variables, $N=2$.
Then $\euf{F}\sb{\lambda,\sigma}$ could be written as
\begin{equation}
\int\sb{\R\times\S}\int\sb{Y}
e\sp{i\phi(x,\tau,\alpha,y)}
\rho(\phi\sb{\alpha\alpha}''(x,1,\alpha,y)/\sigma)
\beta(\tau/\lambda)
a(x,\tau,\alpha,y)
u(y)\,\tau d\tau\,d\alpha\,dy,
\end{equation}
where $a(x,\tau,\alpha,y)$ is a classical symbol
of order $d=\mu+\frac n 2 -1$.
For (\ref{l1-inf-off-caustics}),
we need the bound
\begin{equation}\label{i-i}
\Abs{
\int\sb{\R\times\S}
e\sp{i\phi(x,\tau,\alpha,y)}
\rho(\phi\sb{\alpha\alpha}''(x,1,\alpha,y)/\sigma)
\beta(\tau/\lambda)
a(x,\tau,\alpha,y)\,\tau d\tau\,d\alpha
}
\le C\lambda\sp{\mu+\frac{n+1}{2}}\sigma\sp{-\frac 1 2},
\end{equation}
uniformly in $x$ and $y$ from a small open neighborhood
in $X\times Y$
and for all $\lambda\ge 1$ and $\sigma\le 1$,
$\sigma\ge\lambda\sp{-\frac{m}{m+2}}$.
For simplicity, we assume that $\mu=-n/2$,
so that $\tau a(x,\tau,\alpha,y)$
is a symbol of order zero,
which we denote by $b(x,\tau,\alpha,y)$.
This classical symbol has the development
\begin{equation}
b(x,\tau,\alpha,y)
\sim b\sb 0(x,\alpha,y)+\sum\sb{j\in\N}b\sb j(x,\alpha,y)\tau\sp{-j},
\qquad \tau\ge 1.
\end{equation}
Denote
\begin{equation}\label{i-lambda}
I\sb{\lambda,\sigma}(x,y)=\lambda\sp{-1/2}\int\sb{\R} d\tau
\int\sb{K} d\alpha\,
e\sp{i\phi(x,\tau,\alpha,y)}
b(x,\tau,\alpha,y)
\beta(\tau/\lambda)
\beta(\phi\sb{\alpha\alpha}''(x,1,\alpha,y)/\sigma),
\end{equation}
where $K\subset\S$ denotes $\alpha$-support
of $b(x,\tau,\alpha,y)$.
Substituting $\tau=\lambda z$,
we rewrite $I\sb{\lambda,\sigma}(x,y)$ as
\begin{equation}\label{i-lambda-z}
I\sb{\lambda,\sigma}(x,y)=\lambda\sp{1/2}\int\sb{1/2}\sp{2} dz
\int\sb{K} d\alpha\,
e\sp{i\lambda z \phi(x,1,\alpha,y)}
b(x,\lambda z,\alpha,y)
\beta(z)
\beta(\phi\sb{\alpha\alpha}''(x,1,\alpha,y)/\sigma).
\end{equation}
For (\ref{i-i}), we need the bound
\begin{equation}\label{i-lambda-bound}
\abs{I\sb{\lambda,\sigma}(x,y)}\le C\sigma\sp{-1/2},
\end{equation}
valid for all $\lambda\ge 1$, $\lambda\sp{-m/(m+2)}\le\sigma\le 1$,
and with $C$ independent of $\lambda$ and $\sigma$.

\begin{lemma}\label{lemma-i-bound}
\begin{equation}
\abs{I\sb{\lambda,\sigma}(x,y)}
\le C \sigma\sp{-1/2},
\end{equation}
with $C<\infty$
independent on $\lambda>1$ and $0<\sigma\le 1$.
\end{lemma}

We need this estimate
to be uniform in $\lambda$ and $\sigma$ simultaneously.
Similar estimates were considered in \cite{MR58:563}
and in many other papers.
The result is known to be optimal,
but is not proved in the whole generality
in higher dimensions.
For the sake of completeness, we give our own proof
for the case we are interested in.

\medskip

\begin{proof}
There is a trivial bound
$\abs{I\sb{\lambda,\sigma}(x,y)}\le\lambda\sp{1/2}$,
due to the compact support of the integrand.
This settles the case $0<\sigma\le\lambda\sp{-1}$;
from now on, we assume that $\lambda\sp{-1}\le\sigma\le 1$.

Denote
$b'(x,\tau,\alpha,y)=b(x,\tau,\alpha,y)-b\sb 0(x,\alpha,y)
\in S\sb{\rm cl}\sp{-1}$,
and let
\begin{equation}\label{i-lambda-1}
I\sb{\lambda,\sigma}'(x,y)=\lambda\sp{-1/2}\int\sb{\R} d\tau
\int\sb{K} d\alpha\,
e\sp{i\lambda z\phi(x,1,\alpha,y)}
b'(x,\lambda z,\alpha,y)
\beta(z)
\beta(\phi\sb{\alpha\alpha}''(x,1,\alpha,y)/\sigma).
\end{equation}
Since
$\abs{b(x,\lambda z,\alpha,y)\beta(z)}\le C\lambda\sp{-1}$,
uniformly in $x$, $\alpha$, $\lambda\ge 1$, $1/2\le z\le 2$, and $y$,
there is an easy bound
\[
\abs{I\sb{\lambda,\sigma}'(x,y)}
\le
\lambda\sp{1/2}\int\sb{1/2}\sp{2}
\int\sb{K}
\abs{b'(x,\lambda z,\alpha,y)}
\beta(z)\,dz\,d\alpha
\le C\lambda\sp{-1/2}\le C.
\]
Thus, we only need to consider the bound on (\ref{i-lambda-z})
with $b\sb 0(x,\alpha,y)$ instead of $b(x,\lambda z,\alpha,y)$:
\begin{equation}\label{i-lambda-0}
I\sp{0}\sb{\lambda,\sigma}(x,y)
=\lambda\sp{1/2}\int\sb{\R}
\int\sb{K}
e\sp{i\lambda z\phi(x,1,\alpha,y)}
b\sb 0(x,\lambda z,\alpha,y)
\beta(z)
\beta(\phi\sb{\alpha\alpha}''(x,1,\alpha,y)/\sigma)
\,dz\,d\alpha.
\end{equation}
Denoting by $\hat\beta$ the Fourier transform of $\beta(z)$,
we rewrite $I\sp{0}\sb{\lambda,\sigma}$
as
\begin{equation}
I\sp{0}\sb{\lambda,\sigma}(x,y)
=\lambda\sp{1/2}\int\sb{K}
b\sb 0(x,\alpha,y)
\hat\beta(\lambda\phi(x,1,\alpha,y))
\beta(\phi\sb{\alpha\alpha}''(x,1,\alpha,y)/\sigma)\,d\alpha.
\end{equation}
The statement of the lemma
follows from the bound
\begin{equation}
\abs{I\sp{0}\sb{\lambda,\sigma}(x,y)}\le C\sigma\sp{-1/2},
\end{equation}
which is uniform in $\lambda$, $\lambda\sp{-1}\le\sigma<1$, $x$, and $y$.
We will prove this bound in the next lemma.

This completes the proof of Lemma \ref{lemma-i-bound}.
\end{proof}

\begin{lemma}\label{lemma-sublevel}
Assume that $\phi\sb{\alpha\alpha}''(\alpha)$ vanishes at most of order $m$
on a compact set $K\subset\R$.

If
$f\in L\sp 1(\R)\cap L\sp\infty(\R)$
and
$\beta\in C\sp\infty\sb{0}([1/2,2])$, then
\begin{equation}\label{beta-beta}
\lambda\sp{1/2}\int\sb{K}
f(\lambda\phi(\alpha))\beta(\phi\sb{\alpha\alpha}''(\alpha)/\sigma)\,d\alpha
\le (\norm{f}\sb{L\sp 1}+\norm{f}\sb{L\sp\infty})C\sigma\sp{-1/2},
\end{equation}
uniformly in $\lambda>1$ and $0<\delta\le 1$.
\end{lemma}

Essentially, we are proving the following sublevel
set estimate:
\[
\Abs{\{\alpha\in K\sothat\abs{\phi(\alpha)-\gamma}\le\lambda\sp{-1},
\quad \abs{\phi\sb{\alpha\alpha}''(\alpha)}\ge\sigma\sp{1/2}\}
}
\le \frac{C}{\sqrt{\lambda\sigma}}
\qquad
({\rm uniformly\ in\ }\gamma\in\R).
\]

\begin{proof}
The proof of this estimate is simple, so we can give it in detail.
Let $\rho\in C\sb 0\sp\infty([-2,2])$,
$\rho\at{[-1,1]}\equiv 1$.
We use the partition
\[
1=\rho(\phi\sb\alpha'\sqrt{\lambda/\sigma})
+\left(1-\rho(\phi\sb\alpha'\sqrt{\lambda/\sigma})\right)
\]
to rewrite (\ref{beta-beta}) as a sum of two terms,
\begin{equation}\label{two-terms}
\lambda\sp{\frac 1 2}\int\limits\sb{K}
f(\lambda\phi)
\rho(\phi\sb\alpha'\sqrt{\lambda/\sigma})
\beta(\phi\sb{\alpha\alpha}''/\sigma)\,d\alpha
+
\lambda\sp{\frac 1 2}\int\limits\sb{K}
f(\lambda\phi)
\left(1-\rho(\phi\sb\alpha'\sqrt{\lambda/\sigma})\right)
\beta(\phi\sb{\alpha\alpha}''/\sigma)\,d\alpha,
\end{equation}
which we analyze separately.

The first term in (\ref{two-terms}) is bounded by
\begin{eqnarray}\nonumber
\lambda\sp{\frac 1 2}
\norm{f}\sb{L\sp\infty}
\int\limits\sb{K}
\rho(\phi\sb\alpha'(\alpha)\sqrt{\lambda/\sigma})
\beta(\phi\sb{\alpha\alpha}''(\alpha)/\sigma)\,d\alpha
\\
\le\lambda\sp{\frac 1 2}\norm{f}\sb{L\sp\infty}
\frac{C}
{\sqrt{\lambda/\sigma}\inf\abs{\phi\sb{\alpha\alpha}''}}
\le C\sigma\sp{-\frac 1 2}\norm{f}\sb{L\sp\infty},
\nonumber
\end{eqnarray}
since $\inf\abs{\phi\sb{\alpha\alpha}''}\ge\sigma/2$
on the support of the integrand.
The value of $C$ depends on the bound on the number of
roots of $\phi\sb\alpha'(\alpha)=c$
(this number is bounded uniformly in $c$ due to the
finite type assumption:
$\phi\sb{\alpha\alpha}''$
vanishes of order at most $m$).

The second term is bounded by
\[
\lambda\sp{1/2}\int\sb{K}
f(\lambda\phi(\alpha))
(1-\rho(\phi\sb\alpha'(\alpha)\sqrt{\lambda/\sigma}))
\,d\alpha
\le\lambda\sp{1/2}\norm{f}\sb{L\sp 1}
\frac{C}
{\lambda\inf\abs{\phi\sb{\alpha}'}}
\le C\sigma\sp{-1/2}\norm{f}\sb{L\sp 1},
\]
since $\inf\abs{\phi\sb{\alpha}'}\ge \sqrt{\sigma/\lambda}$
on the support of the integrand.
Again, we need to mention that the
number of roots of
$\phi\sb{\alpha\alpha}''(\alpha)=c$
is bounded uniformly in $c$ due to the finite type assumption.

This proves Lemma \ref{lemma-sublevel}.
\end{proof}

\begin{remark}
The maximal order of vanishing,
$m\in\N$, does not appear in the above lemma.
The statement of the lemma
is also true without the finite type assumption
if we require that $\phi$ is real analytic,
or, more generally,
if we require that $\phi\in C\sp\infty(\R)$
and that $\phi''$ is ``finitely oscillating'' on $K\subset\R$:

\bigskip
{\it
Number of connected components
of the set $\{\alpha\in K\sothat \phi''(\alpha)=c\}$
is bounded uniformly in $c\in\R$.
}
\bigskip

\noindent
This assumption holds for any real analytic function,
but does not hold for all smooth functions;
an example of a smooth function
which is ``infinitely oscillating''
on $[-1,1]$ is $e\sp{-1/x^2}\sin(1/x)$.
\end{remark}
This finishes the proof of (\ref{l1-inf-off-caustics}).

\bigskip

The proof of (\ref{l1-inf-at-caustics}) is similar
but much more straightforward.
One needs to use the following well-known lemma
(see, e.g., \cite{MR2000h:42010}):

\begin{lemma}\label{lemma-vanishing}
If $f(\alpha)$ vanishes at most of order $m$
on $[-2,2]$
and $\beta\in C\sp\infty\sb{0}([-2,2])$,
then
$
\int\sb{\R}
\beta(f(\alpha)/\sigma)\,d\alpha
$
is bounded by $C\sigma\sp{1/m}$.
\end{lemma}

This finishes the proof of
Proposition \ref{prop-l1-inf}.
\end{proof}

\subsection{$L\sp 1 \to L\sp 2$ estimates}

\begin{lemma}
$
\norm{
\euf{F}\sb{\lambda,\sigma}\sp\ast
\euf{F}\sb{\lambda,\sigma}
}\sb{L\sp 1\to L\sp\infty}
+
\norm{
\tilde{\euf{F}}\sb{\lambda,\sigma}\sp\ast
\tilde{\euf{F}}\sb{\lambda,\sigma}
}\sb{L\sp 1\to L\sp\infty}
\le C\lambda\sp{2\mu+n}\sigma\sp{\frac 1 m}.
$
\end{lemma}

\begin{proof}
Since we assume that $\eub C$ is a local graph
(or at least that $\eub C\to Y$ is a submersion,
as in Theorem \ref{theorem-index-m-hardy}),
we can choose the phase function of the form
$\phi(x,\theta,y)=x\cdot\theta-S(\theta,y)$,
with $\theta\in\R^N$, $N=n$,
where $S(\theta,y)$ is homogeneous in $\theta$
of degree $1$.
Then $\theta$ and $y$ can be used as the local coordinates
on $\eub C$.
We can rewrite $\euf{F}$ in the form
\[
\euf{F}\sb{\lambda,\sigma} u(x)
=\int\sb{\R^N\times Y} e\sp{i(x\cdot\theta-S(\theta,y))}a(\theta,y)
\beta(\phi(\theta,y)/\sigma)
u(y)\,d\theta\,dy.
\]
$\euf{F}\sp\ast\euf{F}\in I\sp{-2\mu}(X,X,\varDelta)$
is a Fourier integral operator
with the phase $S(\theta,z)-S(\theta,y)$,
associated to the diagonal
$\varDelta\subset T\sp\ast(Y)\times T\sp\ast(Y)$,
and with $N=n$.

The Fourier integral operator
\begin{eqnarray}
&&\euf{F}\sb{\lambda,\sigma}\sp\ast\euf{F}\sb{\lambda,\sigma} u(z)=
\nonumber
\\
\nonumber
\\
&&\int\limits\sb{\R^N\times Y}
e\sp{i(S(\theta,z)-S(\theta,y))}a(\theta,y)\bar a(\theta,z)
\beta(\mathcal{D}(\theta,y)/\sigma)
\beta(\mathcal{D}(\theta,z)/\sigma)
\beta^2(\abs{\theta}/\lambda)
u(y)\,d\theta\,dy
\nonumber
\end{eqnarray}
has $N=n$ oscillatory variables.
(The number of oscillatory variables cannot be reduced 
since the rank of the matrix
$
\p\sb{\theta\sb i}\p\sb{\theta\sb j}[S(\theta,z)-S(\theta,y)]
$
is zero at $y=z$.)
The order of its symbol is $2\mu$.
This yields the bound
$\const\lambda\sp{2\mu+n}\sigma\sp{\frac 1 m}$
on the $L\sp 1\to L\sp\infty$ action,
with the factor $\sigma\sp{\frac 1 m}$
due to Lemma \ref{lemma-vanishing}.
\end{proof}

\begin{remark}
If $\eub C$ is a local graph, so that
$\det\sb{ij} 
\p\sb{\theta\sb i}\p\sb{y\sb j}S(\theta,y)
\ne 0$,
then $\euf{F}\sp\ast \euf{F}$ is a pseudodifferential operator.
\end{remark}

This lemma yields the following estimate:

\begin{proposition}\label{prop-l1-l2}
Let $\eub C$ and $\euf{F}\in I\sp\mu(X,Y,\eub C)$
be as in Theorem \ref{theorem-index-m-lp-lq},
and $\euf{F}\sb{\lambda,\sigma}$,
$\tilde{\euf{F}}\sb{\lambda,\sigma}$
be given by {\rm (\ref{f-lambda-sigma}), (\ref{f-tilde-lambda-sigma}).}
Then
\begin{equation}
\norm{\euf{F}\sb{\lambda,\sigma}}\sb{L\sp 1\to L\sp{2}}
+\norm{\tilde{\euf{F}}\sb{\lambda,\sigma}}\sb{L\sp 1\to L\sp{2}}
\le C\lambda\sp{\mu+\frac n 2 }\sigma\sp{\frac{1}{2m}}.
\end{equation}

\end{proposition}

\subsection{$L\sp p\to L\sp q$ estimates
for $1<p\le 2\le q <\infty$}

Even when we can not prove the sharp
$\hardy\sp 1\to L\sp q$ estimates
(without the loss of $\epsilon>0$),
we still can prove
the sharp $L\sp p\to L\sp q$ estimates,
for certain values of $p$ and $q$
which also satisfy
$1<p\le 2$, $2\le q <\infty$.
The main tool is the Littlewood-Paley
theory.

We group the pieces $\euf{F}\sb{\lambda,\sigma}$
and $\tilde{\euf{F}}\sb{\lambda,\sigma}$
defined by (\ref{f-lambda-sigma}), (\ref{f-tilde-lambda-sigma})
into $\lambda$-clusters:
\begin{equation}\label{lambda-clusters}
\euf{F}\sb\lambda=\sum\sb{\sigma=2\sp{-j},\,j\in\N,\,\sigma>\sigma\sb 0(\lambda)}
F\sb{\lambda,\sigma}
+\tilde{\euf{F}}\sb{\lambda,\sigma\sb 0(\lambda)}.
\end{equation}
Let us consider this series in the norm of operators
from $L\sp 1 \sb{\mu+\frac n 2+\delta\sb{q}}(Y)$ to $L\sp q(X)$.

\begin{proposition}\label{prop-lp-lq}
Let $\eub C$ and $\euf{F}\in I\sp\mu(X,Y,\eub C)$
be as in Theorem \ref{theorem-index-m-lp-lq},
and let $\euf{F}\sb{\lambda}$
be given by {\rm (\ref{lambda-clusters})}.
Then
\begin{equation}\label{nice-lp-lq}
\norm{\euf{F}\sb{\lambda}}\sb{L\sp 1\to L\sp q}
\le C\lambda\sp{\mu+\frac n 2 +\delta\sb{q}},
\qquad 2\le q<q\sb m,
\end{equation}\label{bad-lp-lq}
\begin{equation}
\norm{\euf{F}\sb\lambda}\sb{L\sp 1\to L\sp q}
\le
C\lambda\sp{\mu+\frac n 2 +\delta\sb{q}
+\kappa\frac{\delta\sb{q}-\delta\sb{q\sb m}}{1/2-\delta\sb{q\sb m}}},
\qquad q\sb m<q<\infty.
\end{equation}
\end{proposition}

\begin{proof}
Interpolating
$L\sp 1\to L\sp\infty$ estimates
from Proposition \ref{prop-l1-inf}
with
$L\sp 1\to L\sp 2$ estimates
from Proposition \ref{prop-l1-l2},
we obtain:

\begin{lemma}
\begin{equation}
\norm{\euf{F}\sb{\lambda,\sigma}}
\sb{L\sp 1\to L\sp q}
\le C\lambda\sp{\mu+\frac n 2 +\delta\sb{q}}
\sigma\sp{\frac 1 {2m}-(\frac 1 m +1)\delta\sb{q}}.
\end{equation}
\end{lemma}
Therefore, (\ref{lambda-clusters})
is dominated by the geometric series.

If $\frac 1 {2m}>(\frac 1 m +1)\delta\sb{q}$
(equivalent with $q<q\sb m$),
then the geometric series is convergent,
and hence is bounded uniformly in $\lambda$.
This proves (\ref{nice-lp-lq}).

If $\frac 1 {2m}<(\frac 1 m +1)\delta\sb{q}$
(equivalent with $q>q\sb m$),
the series (\ref{lambda-clusters}) is dominated by
finitely many terms of the divergent geometric series:
\begin{eqnarray}\label{div-geometric-series}
\norm{\euf{F}\sb{\lambda}}
\sb{L\sp 1\to L\sp q}
&\le&
C\sum\sb{\sigma=2\sp{-j},\,j\in\N,\,\sigma\ge\sigma\sb 0(\lambda)}
\lambda\sp{\mu+\frac n 2 +\delta\sb{q}}
\sigma\sp{-((\frac 1 m +1)\delta\sb{q}-\frac 1 {2m})}
\\
&\le&
C\lambda\sp{\mu+\frac n 2 +\delta\sb{q}}
\sigma\sb 0(\lambda)\sp{-((\frac 1 m +1)\delta\sb{q}-\frac 1 {2m})}.
\nonumber
\end{eqnarray}
Taking into account that
$\sigma\sb 0(\lambda)=\lambda\sp{-\frac{m}{m+2}}$
and that
\begin{eqnarray}
\frac{m}{m+2}
\left(\left(\frac 1 m +1\right)\delta\sb{q}-\frac 1 {2m}\right)
&=&\frac{1}{2(m+2)}
(2(1+m)\delta\sb{q}-1)
\\
&=&\frac \kappa m
(\delta\sb{q}/\delta\sb{q\sb m}-1)
=\kappa\frac{\delta\sb{q}-\delta\sb{q\sb m}}{1/2-\delta\sb{q\sb m}},
\nonumber
\end{eqnarray}
we can rewrite (\ref{div-geometric-series}) in a more convenient form:
\begin{equation}
\norm{\euf{F}\sb{\lambda}}
\sb{L\sp 1\to L\sp q}
\le C
\lambda\sp{\mu+\frac n 2 +\delta\sb{q}
+\kappa\frac{\delta\sb{q}-\delta\sb{q\sb m}}{1/2-\delta\sb{q\sb m}}}.
\end{equation}
This proves (\ref{bad-lp-lq}).
\end{proof}

The estimates stated in Proposition \ref{prop-lp-lq}
can be interpolated with the $L\sp 2\to L\sp 2$-estimates.
If we assume that $\eub C$ is a local graph, then
$\euf{F}\sb\lambda:\;L\sp 2\sb{\mu}\to L\sp 2$,
and we obtain, for $1\le p\le 2$:
\begin{eqnarray}
\norm{\euf{F}\sb\lambda}\sb{L\sp p\to L\sp q}
&\le&
C\lambda\sp{\mu+ n\delta\sb{p}+\delta\sb{q}},
\qquad
2\le q<\frac{2}{1-4\delta\sb{p}\delta\sb{q\sb m}},
\nonumber
\\
\norm{\euf{F}\sb\lambda}\sb{L\sp p\to L\sp q}
&\le&
C\lambda\sp{\mu+n\delta\sb{p}
+(\delta\sb{p}+\delta\sb{q})(1/2+\kappa)
+(\delta\sb{q}-\delta\sb{p})},
\qquad\frac{2}{1-4\delta\sb{p}\delta\sb{q\sb m}}<q<p'.
\nonumber
\end{eqnarray}
According to Littlewood-Paley theory
(\cite{MR94h:35292}, Lemma 2.1),
$\euf{F}=\sum\sb{\lambda=2^l,\,l\in\N}\euf{F}\sb\lambda$
has the same $L\sp p\to L\sp q$ regularity properties
as long as $1<p\le 2\le q<\infty$.
This, together with the duality arguments, proves
Theorem \ref{theorem-index-m-lp-lq}.

\sect{
Microlocal techniques: $\hardy\sp 1\to L\sp q$ estimates
}\label{sect-microlocal-hardy}

We are going to prove
Theorem \ref{theorem-index-m-hardy},
which gives the substitute
of the $L\sp p\to L\sp q$ estimates for $p=1$
($\hardy\sp 1\to L\sp q$ estimates)
and for $q=\infty$
($L\sp p\to \bmo$ estimates).

\subsection{$\hardy\sp 1 \to L\sp\infty$ estimates}

The following is the analogue of
Proposition \ref{prop-l1-inf}.

\begin{proposition}\label{prop-atom-inf}
Let $\eub C\subset T\sp\ast(X)\backslash 0\times T\sp\ast(Y)\backslash 0$
be a smooth canonical relation such that
$\eub C\to X$ is a submersion.
Assume that $\eub C$ has only caustics of the type $A\sb{m+1}$ with $m=1$ or $2$.
Let $\euf{F}\in I\sp\mu(X,Y,\eub C)$ have the polyhomogeneous symbol 
with compact support in $X$, $Y$,
and let $\euf{F}\sb{\lambda,\sigma}$,
$\tilde{\euf{F}}\sb{\lambda,\sigma}$
be given by {\rm (\ref{f-lambda-sigma}), (\ref{f-tilde-lambda-sigma})}.
Then, for any atom $a\sb{\eur Q}$
supported in the cube $\eur Q$ with side $r$, we have
\begin{equation}\label{hardy-inf-off-caustics}
\norm{
\euf{F}\sb{\lambda,\sigma}a\sb{\eur Q}
}\sb{L\sp\infty}
\le C\lambda\sp{\mu+\frac{n+1}{2}}\sigma\sp{-\frac 1 2}
\min(\lambda r,\,(\lambda r)\sp{-1}),
\end{equation}
\begin{equation}
\norm{\euf{F}\sb{\lambda,\sigma}a\sb{\eur Q}}\sb{L\sp\infty}
+
\norm{\tilde{\euf{F}}\sb{\lambda,\sigma}a\sb{\eur Q}}\sb{L\sp\infty}
\le C\lambda\sp{\mu+\frac{n+2}{2}}\sigma\sp{\frac 1 m}
\min(\lambda r,\,(\lambda r)\sp{-1}).
\label{hardy-inf-at-caustics}
\end{equation}
\end{proposition}

\begin{proof}
The proof is similar to \cite{MR92g:35252},
\cite{hardy}.
For the reader's convenience, we reproduce this
proof in Appendix \ref{appendix-hardy}.
We require that $\sigma\ge\lambda\sp{-1/2}$
(equivalent to $m\le 2$)
so that the localizations would not be too fine
and the integration by parts from \cite{MR92g:35252}
could be used verbatim.
\end{proof}

We group the pieces $\euf{F}\sb{\lambda,\sigma}$
into $\lambda$-clusters as in (\ref{lambda-clusters}):
\[
\euf{F}\sb\lambda
=\sum\sb{\sigma=2\sp{-j},\,j\in\N,\,\sigma>\sigma\sb 0(\lambda)}
\euf{F}\sb{\lambda,\sigma}
+\tilde{\euf{F}}\sb{\lambda,\sigma\sb 0(\lambda)}.
\]
The estimates
(\ref{hardy-inf-off-caustics})
and (\ref{hardy-inf-at-caustics})
yield the following bounds on
$\norm{\euf{F}\sb\lambda a\sb{\eur Q}}\sb{L\sp\infty}$:

\begin{corollary}
Assume that $\eub C\to X$ is a submersion
and that $\eub C$ has only caustics of the type $A\sb{m+1}$ with $m=1$ or $2$.
Then
\begin{equation}
\norm{
\euf{F}\sb\lambda a\sb{\eur Q}}\sb{L\sp\infty}
\le C\lambda\sp{\mu+\frac{n+1}{2}+\kappa}\min(\lambda r,\,(\lambda r)\sp{-1}),
\qquad \kappa=\frac 1 2 -\frac{1}{m+2}.
\end{equation}
\end{corollary}

This allows us to conclude that
\begin{equation}
\euf{F}:\; \hardy\sp 1\sb{\mu+\frac{n+1}{2}+\kappa}(Y)
\to L\sp\infty(X).
\end{equation}

\subsection{$\hardy\sp 1 \to L\sp 2$ estimates}

Proposition \ref{prop-atom-inf}
gives the sharp version of
Proposition \ref{prop-l1-inf}.
Now we are going to prove the sharp version
of Proposition \ref{prop-l1-l2}.

\begin{lemma}
Assume that
$\eub C\to Y$ is a submersion
and that $\eub C$ has only caustics of the type $A\sb{m+1}$ with $m=1$ or $2$.
Then, for any atom $a\sb{\eur Q}$
supported in the cube $\eur Q$ with side $r$, we have
\begin{equation}
\norm{
\euf{F}\sb{\lambda,\sigma}\sp\ast
\euf{F}\sb{\lambda,\sigma}a\sb{\eur Q}
}\sb{L\sp\infty}
+\norm{
\tilde{\euf{F}}\sb{\lambda,\sigma}\sp\ast
\tilde{\euf{F}}\sb{\lambda,\sigma}a\sb{\eur Q}
}\sb{L\sp\infty}
\le C\lambda\sp{2\mu+n}\sigma\sp{\frac 1 m}\min(\lambda r,\,(\lambda r)\sp{-1}).
\end{equation}
\end{lemma}

\begin{proof}
The proof is similar
to the proof of Proposition \ref{prop-l1-l2}.
For the $\hardy\sp 1\to L\sp\infty$ estimates,
we can apply the usual machinery
as long as $\min\sigma\approx\lambda\sp{-\frac{m}{m+2}}$
is not smaller than $\lambda\sp{-1/2}$,
that is, as long as $m\le 2$.
\end{proof}

This lemma proves the following sharp version
of Proposition \ref{prop-l1-l2}.

\begin{proposition}\label{prop-atom-l2}
Let $\eub C\subset T\sp\ast(X)\backslash 0\times T\sp\ast(Y)\backslash 0$
be a smooth canonical relation such that
$\eub C\to Y$ is a submersion.
Assume that $\eub C$ has only caustics of the type $A\sb{m+1}$ with $m=1$ or $2$.
Let $\euf{F}\in I\sp\mu(X,Y,\eub C)$ have the polyhomogeneous symbol 
with compact support in $X$, $Y$,
and let $\euf{F}\sb{\lambda,\sigma}$,
$\tilde{\euf{F}}\sb{\lambda,\sigma}$
be given by {\rm (\ref{f-lambda-sigma}), (\ref{f-tilde-lambda-sigma})}.
Then, for any atom $a\sb{\eur Q}$
supported in the cube $\eur Q$ with side $r$, we have:
\begin{equation}\label{hardy-2}
\norm{\euf{F}\sb{\lambda,\sigma}a\sb{\eur Q}}\sb{L\sp{2}}
+\norm{\tilde{\euf{F}}\sb{\lambda,\sigma}a\sb{\eur Q}}\sb{L\sp{2}}
\le C\lambda\sp{\mu+\frac n 2 }\sigma\sp{\frac{1}{2m}}
\min((\lambda r)\sp{1/2},(\lambda r)\sp{-1/2}).
\end{equation}
\end{proposition}

\subsection{$\hardy\sp 1\to L\sp q$ estimates
for small $2\le q<q\sb m$:
$\sigma$-interpolation}

We group the pieces 
$\euf{F}\sb{\lambda,\sigma}$ and $\tilde{\euf{F}}\sb{\lambda,\sigma}$
into $\sigma$-clusters:
\begin{equation}\label{sigma-clusters}
\euf{F}\sb\sigma
=\sum\sb{
\begin{array}{c}
{\scriptstyle \lambda:\;\sigma\ge 2\sigma\sb 0(\lambda)}
\\{\scriptstyle \lambda=2\sp{l},\,l\in\N}
\end{array}
}
\euf{F}\sb{\lambda,\sigma}
+
\sum\sb{
\begin{array}{c}
{\scriptstyle \lambda:\;\sigma\sb 0(\lambda)\le\sigma<2\sigma\sb 0(\lambda)}
\\{\scriptstyle \lambda=2\sp{l},\,l\in\N}
\end{array}
}
\tilde{\euf{F}}\sb{\lambda,\sigma}.
\end{equation}
Then we have
\[
\euf{F}=\sum\sb{\sigma=2\sp{-j},\,j\in\N}\euf{F}\sb\sigma+\euf{F}\sb{\rm nice}.
\]

Proposition \ref{prop-atom-inf}
proves the following bound:

\begin{lemma}
Assume that
$\eub C\to X$ is a submersion
and that $\eub C$ has only caustics of the type $A\sb{m+1}$ with $m=1$ or $2$.
Let $\euf{F}\sb\sigma$ be given by {\rm (\ref{sigma-clusters})}.
Then
\begin{equation}\label{sigma-hardy-inf}
\norm{
\euf{F}\sb\sigma}
\sb{\hardy\sp 1\sb{\mu+\frac{n+1}{2}}\to L\sp\infty}
\le C\sigma\sp{-1/2}.
\end{equation}
\end{lemma}

Proposition \ref{prop-atom-l2}
proves the following:

\begin{lemma}
Assume that
$\eub C\to Y$ is a submersion
and that $\eub C$ has only caustics of the type $A\sb{m+1}$ with $m=1$ or $2$.
Let $\euf{F}\sb\sigma$ be given by {\rm (\ref{sigma-clusters})}.
Then
\begin{equation}\label{sigma-hardy-l2}
\norm{\euf{F}\sb\sigma}
\sb{\hardy\sp 1\sb{\mu+\frac n 2}\to L\sp 2}
\le
C\sigma\sp{\frac 1 {2m}}.
\end{equation}
\end{lemma}

\begin{corollary}
Assume that
both $\eub C\to X$ and $\eub C\to Y$ are submersions
and that $\eub C$ has only caustics of the type $A\sb{m+1}$ with $m=1$ or $2$.
Let $\euf{F}\sb\sigma$ be given by (\ref{sigma-clusters}).
Then
the interpolation of
{\rm (\ref{sigma-hardy-inf})} and {\rm (\ref{sigma-hardy-l2})}
gives
\begin{equation}\label{hardy-q-off-caustics}
\norm{\euf{F}\sb\sigma}
\sb{\hardy\sp 1\sb{\mu+\frac n 2 +\delta\sb{q}}\to L\sp q}
\le
C\sigma\sp{\frac 1 {2m}-(\frac 1 m + 1)\delta\sb{q}},
\qquad 2\le q\le\infty.
\end{equation}
\end{corollary}

The summation $\sum\sb{\sigma=2\sp{-j},\,j\in\N}\euf{F}\sb\sigma$
converges in $\hardy\sp 1\sb{\mu+\frac n 2 +\delta\sb{q}}\to L\sp q$
operator norm (where $q\ge 2$)
if
\[
\frac 1 {2m}>\left(\frac 1 m + 1 \right)\delta\sb{q},
\]
which is equivalent to 
$2\le q<q\sb m$,
$q\sb m=2+\frac 2 m$.
In this case, we conclude that
\begin{equation}\label{no-dependence-on-k}
\euf{F}:\;
\hardy\sp 1\sb{\mu+\frac n 2+\delta\sb{q}}\to L\sp q,
\qquad 2\le q<q\sb m.
\end{equation}
Note that the estimates (\ref{no-dependence-on-k}) do not
depend on the order of caustics.

\subsection{$\hardy\sp 1\to L\sp q$ estimates for $q>q\sb m$:
$\omega$-interpolation}

In the case $m\le 2$, we can derive
the sharp $\hardy\sp 1\to L\sp q$ estimates
for $q>q\sb m$.
According to Proposition \ref{prop-atom-inf},
if $m\le 2$ and
if $a\sb{\eur Q}$ is an atom
supported in the cube $\eur Q$ with side $r$, then
\begin{equation}
\norm{\euf{F}\sb{\lambda,\sigma}a\sb{\eur Q}}\sb{L\sp\infty}
+
\norm{\tilde{\euf{F}}\sb{\lambda,\sigma}a\sb{\eur Q}}\sb{L\sp\infty}
\le C\lambda\sp{\mu+\frac{n+2}{2}}\sigma\sp{\frac 1 m}
\min(\lambda r,\,(\lambda r)\sp{-1}),
\label{hardy-inf-at-caustics-near-2}
\end{equation}
\begin{equation}
\norm{
\euf{F}\sb{\lambda,\sigma}a\sb{\eur Q}
}\sb{L\sp\infty}
\le C\lambda\sp{\mu+\frac{n+1}{2}+\kappa}
\lambda\sp{-\kappa}\sigma\sp{-\frac 1 2}
\min(\lambda r,\,(\lambda r)\sp{-1}).
\end{equation}
According to Corollary \ref{prop-atom-l2},
\begin{equation}
\norm{
\euf{F}\sb{\lambda,\sigma}a\sb{\eur Q}}\sb{L\sp{2}}
+\norm{
\tilde{\euf{F}}\sb{\lambda,\sigma}a\sb{\eur Q}}\sb{L\sp{2}}
\le C\lambda\sp{\mu+\frac n 2-\frac{\kappa}{m}}
\lambda\sp{\frac{\kappa}{m}}\sigma\sp{\frac{1}{2m}}
\min((\lambda r)\sp{1/2},(\lambda r)\sp{-1/2}).
\end{equation}
We introduce a new parameter, $\omega$,
for the values of
$\lambda\sp{-\kappa}\sigma\sp{-\frac 1 2}$
(these values are bounded by $1$
since $\sigma\ge\sigma\sb 0(\lambda)\approx\lambda\sp{-2\kappa}$).
Let us group the operators
$\euf{F}\sb{\lambda,\sigma}$ into $\omega$-clusters
$\euf{F}\sb\omega$, $\omega=2\sp{-k}$, $k\in\N$,
so that
\begin{equation}
\euf{F}=\euf{F}\sb{\rm nice}+\sum\sb{\omega=2\sp{-k},\,k\in\N}
\euf{F}\sb\omega,
\end{equation}
where
\begin{equation}\label{omega-clusters}
\euf{F}\sb\omega
=\sum\sb{
\begin{array}{c}
{\scriptstyle \omega\le
\lambda\sp{-\kappa}\sigma\sp{-1/2} <2\omega}\\
{\scriptstyle\sigma\ge 2\sigma\sb 0(\lambda)}
\end{array}}
{\euf{F}}\sb{\lambda,\sigma}
+\sum\sb{\begin{array}{c}
{\scriptstyle \omega\le
\lambda\sp{-\kappa}\sigma\sp{-1/2} <2\omega}\\
{\scriptstyle\sigma\sb 0(\lambda)\le\sigma< 2\sigma\sb 0(\lambda)}
\end{array}}
\tilde{\euf{F}}\sb{\lambda,\sigma\sb 0(\lambda)},
\end{equation}
$\kappa=\frac 1 2 -\frac{1}{m+2}$,
$\lambda=2^l$, $l\in\N$ and $\sigma=2\sp{-j}$, $j\in\N$.
\begin{lemma}
Assume that
both $\eub C\to X$ and $\eub C\to Y$ are submersions
and that $\eub C$ has only caustics of the type $A\sb{m+1}$ with $m=1$ or $2$.
Let $\euf{F}\sb\omega$ be given by {\rm (\ref{omega-clusters})}.
Then, for any atom $a\sb{\eur Q}$
supported in the cube $\eur Q$ with side $r$, we have
\begin{eqnarray}
&&\norm{
\euf{F}\sb{\omega}}
\sb{\hardy\sp 1\sb{\mu+\frac{n+1}{2}+\kappa} \to L\sp\infty}
\le C\omega,
\\
&&\norm{\euf{F}\sb{\omega}}
\sb{\hardy\sp 1\sb{\mu+\frac n 2-\frac\kappa m}\to L\sp{2}}
\le C\omega\sp{-1/m}.
\label{hardy-2-near-2}
\end{eqnarray}
\end{lemma}

\begin{corollary}
Assume that
both $\eub C\to X$ and $\eub C\to Y$ are submersions
and that $\eub C$ has only caustics of the type $A\sb{m+1}$ with $m=1$ or $2$.
Let $\euf{F}\sb\omega$ be given by {\rm (\ref{omega-clusters})}.
Then
\begin{equation}
\norm{\euf{F}\sb\omega}
\sb{\hardy\sp 1\sb{\mu+\frac n 2
-\frac \kappa m
+\delta\sb{q}\frac{2m+3}{m+2}
}\to L\sp q}
\le
C\omega\sp{-\frac 1 m + 2(\frac 1 m + 1)\delta\sb{q}}.
\end{equation}
\end{corollary}
The series $\sum\sb{\omega=2\sp{-j},\,j\in\N}\euf{F}\sb\omega$
(considered in
$\hardy\sp 1\sb{\mu+\frac n 2-\frac{\kappa}{m}}\to L\sp q$ operator norm)
is dominated
by the geometric series
which is convergent if
$
\frac 1 m<2\left(\frac 1 m + 1 \right)\delta\sb{q},
$
which is equivalent with
$q>q\sb m$.
Therefore,
\begin{equation}
\euf{F}:\;
\hardy\sp 1\sb{\mu+\frac n 2+\delta\sb{q}
+\kappa\frac{\delta\sb{q}-\delta\sb{q\sb m}}{1/2-\delta\sb{q\sb m}}}(Y)
\to L\sp q(X),
\qquad q>q\sb m.
\end{equation}
This finishes the proof of Theorem \ref{theorem-index-m-hardy}.

\sect{Estimates for the half-wave operator}\label{sect-geodesic-flow}

Let $(M,g)$ be a compact Riemann manifold of dimension $n$.
Let $P=\sqrt{-\Delta+1}$, where $\Delta$
is the Laplace operator.
The principal symbol $p(x,\xi)=g\sp{ij}(x)\xi\sb i\xi\sb j$ of $P$
generates the Hamiltonian flow
$\Phi\sb{t}:\;T\sp\ast M\to T\sp\ast M$;
this flow leaves invariant the cosphere bundle
\[
S\sp\ast M=\{(x,\xi)\in T\sp\ast M\sothat p(x,\xi)=1\}.
\]
The geodesics of unit speed on $M$
are the curves $t\mapsto\pi\Phi\sb t((x,\xi))$,
$(x,\xi)\in S\sp\ast M$.
Let
$\pi$ be the canonical projection
$T\sp\ast M\to M$.
We say that the time $t$ is non-conjugate if
the bicharacteristics which start at the moment $t=0$ 
at any point $x\in M$
do not form caustics in time $t$,
so that
$\pi\Phi\sb{t}:\;S\sp\ast M\to M$
is of maximal rank:
\begin{equation}
\rank
d\left(\pi\Phi\sb{t}\at{S\sp\ast\sb{x}M}\right)(\xi)=n-1.
\end{equation}
Here $\xi$ is a point in the fiber
$S\sp\ast\sb{x}M$ of the cosphere bundle at the point $x$.

Assume that at $t=T$ the map
$\pi\Phi\sb{t}:\;S\sp\ast M\to M$
is no longer of maximal rank
at the point $(x,\xi)$, where $\xi\in S\sp\ast\sb x M$:
\begin{equation}
\rank d\left(\pi\Phi\sb{t}\at{S\sp\ast\sb{x}M}\right)(\xi)
<n-1,\qquad t=T.
\end{equation}
The integral kernel $K\sp{t}$ 
of the half-wave operator $e\sp{i t P}$
can be represented
as a finite sum of oscillatory integrals
of the form
\[
K\sp{t}(x,y)=\int\sb{\R^n} e\sp{i(x\cdot \xi-\phi(t,y,\xi))}
a\sb{t}(y,\xi)\,d\xi,
\]
where $a\sb{t}(y,\xi)$
is a classical symbol of order $0$.
(See \cite{MR94c:35178}, Section 4.)
This representation is
valid for $(x,\xi,y)$ supported in a small open conic
neighborhood of $M\times\R^n\times M$
and $t$ in a small open neighborhood of $T$.
We apply our results on $L\sp p\to L\sp q$ estimates
(Theorem \ref{theorem-index-m-lp-lq})
to the half-wave operator
$e^{it P}$
with the integral kernel $K\sp{t}(x,y)$.

\begin{theorem}\label{theorem-lp-lq-half-wave}
If for $t\le T$ the geodesic flow $\Phi\sb{t}$
forms only caustics of the type $A\sb{m'+1}$ with $m'\le m$,
then for $0<t\le T$
and for $1<p\le p<\infty$
such that $(p,q)\sp\dagger\notin\overline{\euf{C}\sb m}$
the $L\sp p\to L\sp q$ estimates are caustics-insensitive.

Precisely,
\begin{eqnarray}
&&
e\sp{i t P}P\sp{-n\delta\sb{p}-\delta\sb{q}}:\;
L\sp p\to L\sp q,
\qquad (p,q)\sp\dagger\in\euf{A}\sb m,
\label{good-1-1}
\\
&&
e\sp{i t P}P\sp{-n\delta\sb{q}-\delta\sb{p}}:\;
L\sp p\to L\sp q,
\qquad (p,q)\sp\dagger\in\euf{B}\sb m.
\label{good-inf-inf}
\end{eqnarray}

\noindent
For $(p,q)\sp\dagger\in\euf{C}\sb m$,
the estimates depend on the order of the caustic,
which is given by $\kappa=\frac 1 2 -\frac{1}{m+2}$:
\begin{eqnarray}
&&e\sp{i t P}P\sp{
-n\delta\sb{p}
-(\delta\sb{p}+\delta\sb{q})(1/2+\kappa)
-(\delta\sb{q}-\delta\sb{p})
}:\;
L\sp p\to L\sp q,
\quad
(p,q)\sp\dagger\in\euf{C}\sb m,
\ q\le p',
\label{bad-q-le-pp}
\\
&&e\sp{i t P}P\sp{
-n\delta\sb{q}
-(\delta\sb{p}+\delta\sb{q})(1/2+\kappa)
-(\delta\sb{p}-\delta\sb{q})
}:\;L\sp p\to L\sp q,
\quad
(p,q)\sp\dagger\in\euf{C}\sb m,
\ q>p'.
\label{bad-q-ge-pp}
\end{eqnarray}
The regions $\euf{A}\sb m$, $\euf{B}\sb m$,
and $\euf{C}\sb m$ in $(1/p,1/q)$-plane
are defined in Definition \ref{def-omega-regions} 
(see also Figure \ref{fig-omega-regions}).
\end{theorem}

We can use these results
to investigate precisely 
the blow-up of the solution
just before the formation of the caustics.
At non-conjugate times $t$,
the estimates on the half-wave operator
$e^{i t P}$ 
are given by the estimates 
(\ref{bad-q-le-pp}), (\ref{bad-q-ge-pp})
with $\kappa=0$.
As $t$ approaches the moment $T$ when
the geodesic flow starts forming caustics,
these estimates blow up
(and the estimates with nonzero $\kappa$
are to be used).
As was shown in \cite{MR1829351},
if $T>0$ is such that $t$ is non-conjugate
for $t\in(T-\epsilon,T)$, for some $\epsilon>0$,
then the $L\sp{q'}\to L\sp q$ estimates on
the half-wave operator $e\sp{i t P}$ may blow up as $t\to T$
at most as
\begin{equation}
\norm{e\sp{i t P}P\sp{-(n+1)/2}}
\sb{L\sp{q'}\to L\sp{q}}
\le C\sb{q,M}(T)\abs{T-t}\sp{-K(n-1)\delta\sb{q}},
\qquad 2\le q<\infty,
\end{equation}
where $K=4$.
This is an a priori value; $K$
could be shown to be smaller when the geodesic flow
forms some particular caustics.

\begin{theorem}\label{theorem-lp-lpp-blow-up}
Let $T>0$ and suppose there exists $\epsilon>0$
such that $t$ is non-conjugate for
$T-\epsilon\le t<T$.
Assume that for $T\le t\le T+\epsilon$ the geodesic flow $\Phi\sb{t}$
forms only simple caustics of index at most $m$
(e.g. caustics of the type $A\sb{m+1}$).
Let $2\le q<\infty$, $1/q+1/q'=1$.
We have for $T-\epsilon/2\le t<T$:
\begin{eqnarray}
&&
\norm{e\sp{i t P}P\sp{-(n+1)\delta\sb{q}}}
\sb{L\sp{q'}(M)\to L\sp q(M)}\le C\sb{q,M}(T)\abs{T-t}\sp{-\delta\sb{q}},
\qquad 2\le q<\infty,
\\
&&\norm{e\sp{i t P}P\sp{-(n+1)\delta\sb{q}-2\kappa\delta\sb{q}}}
\sb{L\sp{q'}(M)\to L\sp q(M)}\le C\sb{q,M}(T),
\qquad 2\le q<\infty,
\end{eqnarray}
where
$\delta\sb{q}=\frac 1 2-\frac 1 q$
and
$\kappa=\frac 1 2 -\frac 1 {m+2}$.
\end{theorem}
\begin{proof}
We reduce the number of oscillatory variables 
in the representation of $K\sp{t}$ to $2$, 
which is possible in an open neighborhood of simple caustics,
and use the polar coordinates
$(\lambda,\alpha)\in\R\sb{+}\times \S$
in the $\theta$-space.
We exploit the fact that
$\abs{\phi\sb{\alpha\alpha}''}\ge \const\abs{T-t}$
if $t$ is non-conjugate for $T-\epsilon<t<T$.
(This bound is easy for stable caustics.
For the generic situation, see \cite{MR1829351}, Lemma 2.4.)
We also use the Littlewood-Paley decomposition
for $K\sp{t}$
(to interpolate
$L\sp 1\to L\sp\infty$ estimates
on $K\sp{t}\sb\lambda$
with $L\sp 2\to L\sp 2$ estimates).
The rest of the theorem 
is the same as the proof of the statement (\ref{l1-inf-off-caustics})
of Proposition \ref{prop-l1-inf}.
Again, the optimal estimate
with the factor $\abs{\det\phi\sb{\alpha\alpha}''}^{-1/2}$
for the oscillatory integral
is readily available since $\alpha$ is one-dimensional.
\end{proof}

The interpolation of the $L\sp p\to L\sp q$ estimates 
which remain valid at the caustics
(Theorem \ref{theorem-lp-lq-half-wave})
and the asymptotics
which describe the blow-up of the usual $L\sp p\to L\sp {p'}$
estimates (Theorem \ref{theorem-lp-lpp-blow-up})
gives the complete description of the
behavior of the blow-up of
$L\sp p\to L\sp q$ estimates
just before the geodesic flow forms caustics.

\appendix 

\sect{Consistency of the definition of $\mathcal{D}$}\label{appendix-consistency}

In this section we prove Lemma \ref{lemma-consistency}:

\bigskip

\noindent
\begin{lemma}[Lemma \ref{lemma-consistency}]
Let $\Lambda$ be a smooth closed conic Lagrangian submanifold
of $T\sp\ast(X)\backslash 0$.
Let $\phi(x,\theta)\in C\sp\infty(X\times\R^N)$
be a smooth non-degenerate phase function
which parametrizes $\Lambda$:
\[
\Lambda=\{(x,d\sb{x}\phi(x,\theta))
\sothat d\sb\theta(x,\theta)=0\}.
\]
Let $\alpha=\{\alpha\sb i\}$, $1\le i\le N-1$,
be local coordinates on the unit sphere $\S\sp{N-1}$.
We use $(\lambda,\alpha)\in \R\sb{+}\times\S\sp{N-1}$
as local coordinates in $\R^N$.
Then
$\mathcal{D}
=\det\sb{ij}(\lambda\sp{-1}\phi\sb{\alpha\sb i\alpha\sb j}''\at{\Lambda})$,
$1\le i,\,j\le N-1$,
is a smooth function on $\Lambda$ defined up to a nonzero factor:
\[
\mathcal{D}\in
C\sp\infty(\Lambda)\slash C\sp\infty\sb\times(\Lambda).
\]

\end{lemma}

We split the proof into two parts:
In the first part, we will show that
if we use the maximal number of oscillatory
variables,
then $\mathcal{D}$
is defined up to a nonzero factor.
In the second part, we show that
$\mathcal{D}$ is multiplied by a nonzero factor
if we reduce the number of oscillatory variables.

\medskip

($\imath$)\quad 
Let us check that, up to a factor, $\mathcal{D}$
does not depend on the chosen parametrization
of $\Lambda$
if we use the maximal number $N=n$ of oscillatory
variables.
$\Lambda$
can be parametrized (locally) by $\theta$
with $\theta\in\R^N$, $N=n$.
Assume there are two different phase functions
$\phi(x,\theta)$
and $\psi(x,\vartheta)$,
$\theta\in\R^N$, $\vartheta\in\R^N$,
and that
both $\theta$ and $\vartheta$
can be used as local coordinates on $\Lambda$.
According to e.g. \cite{MR96m:58245},
there is a function $g(x,\theta)$, homogeneous of
degree $1$ in $\theta$, such that
$\phi(x,\theta)=\psi(x,g(x,\theta))$.
We rewrite $\phi$ and $\psi$ as
\begin{equation}
\phi=\phi(x,\lambda,\alpha),
\qquad
\psi=\psi(x,\tau,\beta),
\end{equation}
where $\lambda=\abs{\theta}$, $\tau=\abs{\vartheta}$,
and $\alpha$, $\beta$ are local coordinates on $\mathbb{S}\sp{N-1}$.
Then there is a smooth function $\beta(x,\alpha)$
and a smooth function $c(x,\alpha)\ne 0$
such that
\begin{equation}\label{phi-vs-psi}
\phi(x,\lambda,\alpha)=\psi(x,c(x,\alpha)\lambda,\beta(x,\alpha)).
\end{equation}
To simplify the notations, 
we will assume the summation with respect to the repeating indices 
and will not write the subscripts of $\alpha$, $\beta$, \dots\ at all, 
assuming that e.g. 
$\det A\sb{\alpha\alpha'}''$ stands for $\det\sb{ij} \p\sb{\alpha\sb i}\p\sb{\alpha\sb j}A$
and 
$A\sb{\alpha\alpha'}''\,d\alpha\,d\alpha'$
stands for
$\sum\sb{ij}\p\sb{\alpha\sb i}\p\sb{\alpha\sb j}A\,d\alpha\sb i\,d\alpha\sb j$.

We differentiate
relation (\ref{phi-vs-psi}) twice with respect to $\alpha$:
\[
\phi\sb{\alpha\alpha'}''
=\psi\sb{\beta\beta'}''J\sp\beta\sb\alpha J\sp{\beta'}\sb{\alpha'}
+\psi\sb{\beta}'J\sp\beta\sb{\alpha\alpha'}
+\psi\sb{\beta\lambda}''J\sp\beta\sb{\alpha}c\sb{\alpha'}
+\psi\sb{\beta\lambda}''J\sp\beta\sb{\alpha'}c\sb{\alpha}
+\psi\sb{\lambda\lambda}''c\sb{\alpha}^2
+\psi\sb\lambda'c\sb{\alpha\alpha},
\]
where $J\sp\beta\sb\alpha(x,\alpha)=\p\beta(x,\alpha)/\p\alpha$
(that is,
$J\sp{\beta\sb j}\sb{\alpha\sb i}(x,\alpha)
=\p{\beta\sb j}(x,\alpha)/\p{\alpha\sb i}$).
Taking into account that $\psi\sb{\lambda\lambda}''\equiv 0$
($\psi$ is homogeneous of degree $1$ in $\lambda$),
while $\psi\sb\lambda'$, $\psi\sb{\beta}'$,
and
$\psi\sb{\lambda\beta}''(x,\lambda,\beta)
=\lambda\sp{-1}\psi\sb{\beta}'(x,\lambda,\beta)$
vanish identically on $\Lambda$,
we deduce that
\[
\det\phi\sb{\alpha\alpha'}''(x,\lambda,\alpha)
=(\det J\sp\beta\sb{\alpha}(x,\alpha))^2
\det\psi\sb{\beta\beta'}''(x,c(x,\alpha)\lambda,\beta(x,\alpha)),
\]
where 
$\det J\sp\beta\sb{\alpha}(x,\alpha)
=\det \p\sb\alpha\beta(x,\alpha)
=\det\sb{ij} \p\sb{\alpha\sb i}\beta\sb j(x,\alpha)
\ne 0$.

\medskip

($\imath\imath$)\quad 
Let us check that
$\mathcal{D}$
as an element of $C\sp\infty(\Lambda)/C\sp\infty\sb\times(\Lambda)$
is not affected by the reduction of oscillatory variables.
We consider the phase function $\phi(x,\lambda,\alpha)$.
Assume that the coordinates $\alpha$
split into $\alpha=(\rho,\sigma)$
so that
$\phi\sb{\sigma\sigma'}''$ is non-degenerate.
Then there exists a smooth function
$\varSigma(x,\rho)$ such that
the condition $\phi\sb\sigma'(x,\lambda,\rho,\sigma)=0$
is equivalent with $\sigma=\varSigma(x,\rho)$.
The phase function
\begin{equation}\label{phase-psi}
\psi(x,\lambda,\rho)
=\phi(x,\lambda,\rho,\varSigma(x,\rho))
\end{equation}
parametrizes the same canonical relation as $\phi$ does.
We are going to prove that $\det\phi\sb{\alpha\alpha'}''$
and $\det\psi\sb{\rho\rho'}''$ differ
by a nonzero factor (namely, $\det\phi\sb{\sigma\sigma'}''$).

In what follows, we drop off the dependence
on $x$ and $\lambda$.
Differentiating (\ref{phase-psi}) with respect to $\rho$, we get
\begin{equation}
\psi\sb{\rho}'(\rho)
=\phi\sb{\rho}'(\rho,\varSigma(\rho))
+\phi\sb\sigma'(\rho,\varSigma(\rho))
J\sp{\sigma}\sb{\rho}(\rho),
\end{equation}
where $J\sp{\sigma}\sb{\rho}(\rho)=\p\varSigma(\rho)/\p\rho$.
\begin{eqnarray}\label{psi-alpha-alpha}
\psi\sb{\rho\rho'}''(\rho)
&=&\phi\sb{\rho\rho'}''(\rho,\varSigma(\rho))
+\phi\sb{\rho\sigma}''(\rho,\varSigma(\rho))
g\sp{\sigma}\sb{\rho'}(\rho)
+\phi\sb{\rho'\sigma}''(\rho,\varSigma(\rho))
J\sp{\sigma}\sb{\rho}(\rho)
\\
&&
+
\phi\sb{\sigma\sigma'}''(\rho,\varSigma(\rho))
J\sp{\sigma}\sb{\rho}(\rho)
J\sp{\sigma'}\sb{\rho'}(\rho)
+\phi\sb\sigma'(\rho,\varSigma(\rho))
J\sp{\sigma}\sb{\rho\rho'}(\rho).
\nonumber
\end{eqnarray}
The last term in the right-hand side of (\ref{psi-alpha-alpha})
vanishes identically on the canonical relation
(where $\phi\sb\alpha'=(\phi\sb\rho',\phi\sb\sigma')\equiv 0$).
Using the identity
\[
0\equiv
\p\sb{\rho}(\phi\sb\sigma'(\rho,\varSigma(\rho)))
=
\phi\sb{\rho\sigma}''(\rho,\varSigma(\rho))
+\phi\sb{\sigma\sigma'}''(\rho,\varSigma(\rho))
J\sp{\sigma'}\sb{\rho}(\rho),
\]
we can express
$J\sp{\sigma}\sb{\rho}(\rho)
=-\phi\sb{\rho\sigma'}''(\rho,\varSigma(\rho))\phi\sp{\sigma\sigma'}(\rho,\varSigma(\rho))$,
where $\phi\sp{\sigma\sigma'}(\rho,\sigma)$
denotes the matrix inverse to $\phi\sb{\sigma\sigma'}''(\rho,\sigma)$.
We rewrite (\ref{psi-alpha-alpha})
as
\begin{equation}\label{psi-alpha-alpha-2}
\psi\sb{\rho\rho'}''(\rho)
=\phi\sb{\rho\rho'}''(\rho,\varSigma(\rho))
-
\phi\sb{\sigma\rho}''(\rho,\varSigma(\rho))
\phi\sp{\sigma\sigma'}(\rho,\varSigma(\rho))
\phi\sb{\sigma'\rho'}''(\rho,\varSigma(\rho)).
\end{equation}
To compute the determinant of (\ref{psi-alpha-alpha-2}),
we use the identity
\begin{equation}\label{determinant-identity}
\det(A-B D\sp{-1}C)\det D
=\det{\left[\begin{array}{cc}A&B\\C&D\end{array}\right]},
\end{equation}
where
$A$ and $D$ are square matrices and $\det D\ne 0$,
which follows from the matrix identity
\[
\left[\begin{array}{cc}A-B D\sp{-1}C&0\\0&D\end{array}\right]
=\left[\begin{array}{cc}I&-B D\sp{-1}\\0&I\end{array}\right]
\left[\begin{array}{cc}A&B\\C&D\end{array}\right]
\left[\begin{array}{cc}I&0\\-D\sp{-1}C&I\end{array}\right].
\]
Identity (\ref{determinant-identity})
allows to write the determinant of (\ref{psi-alpha-alpha-2})
in the form of the desired relation:
\begin{equation}
\det\psi\sb{\rho\rho'}''\det\phi\sb{\sigma\sigma'}''
=\det\phi\sb{\alpha\alpha'}'',
\qquad {\rm where}\quad {\det\phi\sb{\sigma\sigma'}''}\ne 0.
\end{equation}

This finishes the proof of Lemma \ref{lemma-consistency}.

\sect{$\hardy\sp 1\to L\sp\infty$
estimates on $(\lambda,\sigma)$-pieces}\label{appendix-hardy}
\label{sect-sharp-estimates}

In this section we prove the following Lemma,
needed for 
sharp $\hardy\sp 1\to L\sp q$ 
and $L\sp p\to \bmo$ estimates.

\begin{lemma}\label{key-lemma-hardy-l-infinity}
Let $\euf{F}\in I\sp\mu(X,Y,\eub C)$
be associated to a canonical relation
such that
$\eub C\to X$ is a submersion,
and let
$\eub C$ have only caustics of the type $A\sb{m+1}$ with $m=1$ or $2$.
Then, for any atom $a\sb{\eur Q}$
supported in the cube $\eur Q$ with side $r$, we have
\begin{equation}\label{hardy-to-l-infinity-estimates-on-f-pieces}
\norm{\Lambda\sp{-\mu-(n+N)/2}\euf{F}\sb{\lambda,\sigma}
a\sb{\eur Q}}\sb{L\sp\infty(X)}
\le\const\sigma\sp{1/m}
\min(\lambda r,(\lambda r)\sp{-1}).
\end{equation}

\end{lemma}

For simplicity, we consider
$\euf{F}\in I\sp{\mu}(X,Y,\eub C)$
with $\mu=-(N+n)/2$.
This implies that $a(x,\theta,y)\in S^d$ with
$d=-(n+N)/2-(N-n)/2=-N$.
Let $a\sb{\eur Q}$ be an atom
supported in the cube $\eur Q$ with side $r$
(following \cite{MR92g:35252}, we may assume that $r\le 1$).
We want to show that
for any $x$,
\begin{equation}
\abs{\euf{F}\sb{\lambda,\sigma} a\sb{\eur Q}(x)}
\le\const\sigma\sp{1/m}\min(\lambda r,(\lambda r)\sp{-1}).
\end{equation}

We will decompose and bound the pieces $\euf{F}\sb{\lambda,\sigma}$
following the discussion on pages 238-241
in \cite{MR92g:35252}.
For a particular $\lambda$,
we introduce unit vectors $\theta\sb\lambda\sp{\nu}$,
with
$1\le\nu\le N(\lambda\sp{-1/2})\approx\lambda\sp{\frac{m-1}{2}}$,
equidistributed on the unit sphere in the $\theta$-space $\R^N$,
so that
$\abs{\theta\sb\lambda\sp{\nu}
-\theta\sb\lambda\sp{\nu\sp\prime}}\ge\const\lambda\sp{-\frac{1}{2}}$
for $\nu\ne\nu\sp\prime$.
We introduce a corresponding partition
of unity,
\[
1=\sum\sb{\nu=1}\sp{N(\lambda\sp{-1/2})}
\psi\sb\lambda\sp{\nu}(\theta),
\]
where the functions
$\psi\sb\lambda\sp{\nu}$
are homogeneous of degree $0$
and supported in the spherical angles
$\Omega\sb\lambda\sp{\nu}$
with the span $\sim\lambda\sp{-1/2}$,
centered at $\theta\sb\lambda\sp{\nu}$:
\[
\psi\sb\lambda\sp{\nu}(\theta)\ne 0
\quad{\rm only\ if}\quad
\Abs{\frac{\theta}{\abs{\theta}}-\theta\sb\lambda\sp{\nu}}
\le\const\lambda\sp{-\frac{1}{2}}.
\]
We assume that
$
\abs{\p\sb{\theta}\sp{\alpha}
\psi\sb\lambda\sp{\nu}(\theta)}
\le \const\lambda\sp{\frac{\abs{\alpha}}{2}}
\abs{\theta}\sp{-\abs{\alpha}}.
$

We denote the integral kernels of
$\euf{F}\sb{\lambda,\sigma}$,
$\,\tilde{\euf{F}}\sb{\lambda,\sigma}$
by $K\sb{\lambda,\sigma}(x,\theta,y)$ and
$\tilde{K}\sb{\lambda,\sigma}(x,\theta,y)$.
We introduce $\euf{F}\sb{\lambda,\sigma}\sp{\nu}$
by
\[
\euf{F}\sb{\lambda,\sigma}\sp{\nu}u(x)
=\int
K\sb{\lambda,\sigma}\sp{\nu}(x,\theta,y)
u(y)\,d\theta\,dy,
\]
where
$
K\sb{\lambda,\sigma}\sp{\nu}(x,\theta,y)
=
\psi\sb\lambda\sp{\nu}(\theta)
K\sb{\lambda,\sigma}(x,\theta,y).
$

  From now on, we assume that
$x\in X$ is fixed.
We need to introduce the ``exceptional set''
associated to $x$.
According to \cite{hardy},
the assumption that
$\eub C\to X$ is a submersion
allows one to choose the phase $\phi$
in the form
\[
\phi(x,\theta,y)
=\langle
G(x,\theta,y\sp{\prime\prime})-y\sp\prime,
\theta\rangle,
\]
where
$y=(y\sp\prime,y\sp{\prime\prime})\in\R\sp{N}\times\R\sp{n-N}=Y$
are certain local coordinates.
For given $\lambda$ and $\nu$, we define
$\mathcal{R}\sp{\nu}\sb{x,\lambda\sp{-1}}\subset Y$
by
\begin{equation}
\mathcal{R}\sp{\nu}\sb{x,\lambda\sp{-1}}
=\{y \sothat\,
\abs{\langle
G(x,\theta\sb\lambda\sp\nu,y\sp{\prime\prime})-y\sp\prime,
\theta\sb\lambda\sp{\nu}
\rangle}
\le \lambda\sp{-1},\,
\abs{G(x,\theta\sb\lambda\sp\nu,y\sp{\prime\prime})
-y\sp\prime}
\le \lambda\sp{-\frac{1}{2}}
\},
\end{equation}
with
$\abs{\mathcal{R}\sp{\nu}\sb{x,\lambda\sp{-1}}}
\le\const\lambda\sp{-1}\cdot\lambda\sp{-\frac{m-1}{2}}$.
We set
$\chi\sb{\mathcal{R}\sp{\nu}\sb{x,\lambda\sp{-1}}}(x)
$
to be the characteristic function of
$\mathcal{R}\sp{\nu}\sb{x,\lambda\sp{-1}}$.

Let $a\sb{\eur Q}$ be an atom
supported in the cube $\eur Q$ with side $r$:
\[
\abs{\eur Q}=r^n,
\qquad
\norm{a\sb{\eur Q}}\sb{L\sp\infty}\le r\sp{-n},
\qquad
\norm{a\sb{\eur Q}}\sb{L\sp 1}\le 1,
\qquad
\int\sb{\eur Q} a\sb{\eur Q}=0.
\]
We consider
\begin{eqnarray}
\euf{F}\sb{\lambda,\sigma}a\sb{\eur Q}(x)
&=&
\sum\sb{\nu}
\int
\chi\sb{\mathcal{R}\sp{\nu}\sb{x,\lambda\sp{-1}}}(y)
K\sb{\lambda,\sigma}\sp{\nu}(x,\theta,y)
a\sb{\eur Q}(y)\,d\theta\,dy
\nonumber
\\
&+&
\sum\sb{\nu}
\int
(1-\chi\sb{\mathcal{R}\sp{\nu}\sb{x,\lambda\sp{-1}}}(y))
K\sb{\lambda,\sigma}\sp{\nu}(x,\theta,y)
a\sb{\eur Q}(y)\,d\theta\,dy.
\label{3.4}
\end{eqnarray}
We need to know the absolute value of this expression.

({\it i})
The absolute value of the first term in the right-hand side
is bounded by
\begin{eqnarray}
&&\sum\sb{\nu}\int
\chi\sb{\mathcal{R}\sp{\nu}\sb{x,\lambda\sp{-1}}}(y)
\Abs{K\sb{\lambda,\sigma}\sp{\nu}(x,\theta,y)a\sb{\eur Q}(y)}\,d\theta\,dy
\nonumber
\\
&&\le
C\lambda\sp{-N}
\sum\sb{\nu}
\int\chi\sb{\mathcal{R}\sp{\nu}\sb{x,\lambda\sp{-1}}}(y)
\psi\sp{\nu}\sb\lambda(\theta)
\beta({\abs{\theta}}/{\lambda})
\beta(\mathcal{D}(x\sp{\prime\prime},\theta,y)/\sigma)
\abs{a\sb{\eur Q}(y)}\,d\theta\,dy.
\label{3.5}
\end{eqnarray}

$\circ$
In (\ref{3.5}), we have already applied the bound
$C\lambda\sp{-N}$
on the symbol $a(x,\theta,y)\in S\sp{-N}$
at $\abs{\theta}\sim\lambda$.

$\circ$
Summation in $\nu$ converges since
$\sum\sb{\nu}\psi\sb\lambda\sp{\nu}(\theta)=1$.

$\circ$
If $d\sb\theta\mathcal{D}\ne 0$ ($m=1$),
then the integration in $\theta$
contributes $\const\sigma\lambda^N$,
where
$\sigma$ appears due to the support properties
of $\beta(\mathcal{D}/\sigma)$.

More generally,
assume that at a point $p\in \eub C$
there is a simple caustic of the type $A\sb{m+1}$ with $m=1$ or $2$.
Then there is a vector field $V=a\sb{j}\p\sb{\theta\sb j}$,
$V\in C\sp\infty(\Gamma(T(\eub C)))$,
such that
$V^m \mathcal{D}\at{p}\ne 0$.
We define $\Theta=\frac{\theta}{\lambda}\in\R^N$,
so that the region of integration
in $\Theta$ is bounded uniformly in $\lambda$.
Note that $d\theta=\lambda\sp{N}\,d\Theta$.
We can choose the coordinates
so that
$\p\sb{\theta\sb N}^m \mathcal{D}\ne 0$,
in an open neighborhood
of $p$.
The expression $\p\sb{\Theta\sb N}^m\mathcal{D}$
is homogeneous of degree zero in $\lambda$,
so that
$\abs{\p\sb{\Theta\sb N}^m\mathcal{D}}\ge\const>0$
uniformly in $\lambda$, $\sigma$.
Therefore, 
\[
\int\sb{\R}
\beta(\mathcal{D}/\sigma)\,d\Theta\sb N
\le\const\sigma\sp{\frac{1}{m}}.
\]
The integration in $\Theta\sb{1},\dots,\Theta\sb{N-1}$
converges since the support of (\ref{3.5})
in
$\Theta=\frac{\theta}{\lambda}$
is bounded (uniformly in $\lambda$, $\sigma$).
We conclude that
the integration in $\theta$
contributes $\const\lambda^N\sigma\sp{\frac{1}{m}}$.

$\circ$
Finally,
due to the bound
$\norm{a\sb{\eur Q}}\sb{L\sp\infty(Y)}
\le\abs{\eur Q}\sp{-1}$
together with the support properties of
$a\sb{\eur Q}$
and
$\chi\sb{\mathcal{R}\sp{\nu}\sb{x,\lambda\sp{-1}}}(y)$,
the integral in $y$
contributes the factor
$\min(1,(\lambda r)\sp{-1})$.

Taking the product of all of the above factors, we
obtain
$\const\sigma\sp{\frac{1}{m}}\min(1,(\lambda r)\sp{-1})$.

({\it ii})
For the absolute value of the second term
in the right-hand side of $(\ref{3.4})$ we have:

$\circ$
In each $\nu$-term,
we can integrate by parts as in \cite{MR92g:35252}
(we need the assumption $\sigma\ge\lambda\sp{-\frac{1}{2}}$
to obtain an analogue of the inequalities $(3.19)$
in \cite{MR92g:35252};
the argument is the same as theirs),
getting the factor
\[
\left(
1+\lambda^2\abs{\langle
G(x,\theta\sb\lambda\sp{\nu},y\sp{\prime\prime})
-y\sp\prime,
\theta\sb\lambda\sp{\nu}
\rangle
}^2
+\lambda\abs{
G(x,\theta\sb\lambda\sp{\nu},y\sp{\prime\prime})
-y\sp\prime
}^2
\right)\sp{-M},
\qquad {\rm for\ any}\quad M\in\N.
\]
The integral of the product
of this expression with $a\sb{\eur Q}(y)$
with respect to $y$
contributes the same factor
$\min(1,(\lambda r)\sp{-1})$
as above.
The rest of the analysis is the same as
for the first term in the right-hand side of $(\ref{3.4})$.

We conclude that
$\abs{\euf{F}\sb{\lambda,\sigma}a\sb{\eur Q}(x)}
\le\const\sigma\sp{\frac{1}{m}}\min(1,(\lambda r)\sp{-1})$.

The bound in the case
$\lambda r<1$ follows from \cite{hardy}.
Let us recall the argument.
We fix some point $\bar y\in \eur Q$.
Since $\int a\sb{\eur Q}(y)dy=0$, we can write
\begin{eqnarray}\label{extra-lambda-r}
\euf{F}\sb{\lambda,\sigma} a\sb{\eur Q}(x)
&=&\int
\left[
K\sb{\lambda,\sigma}(x,\theta,y)
-K\sb{\lambda,\sigma}(x,\theta,\bar y)
\right]
a\sb{\eur Q}(y)d\theta\,dy
\nonumber \\
&=&
\int\sb{0}\sp 1dt\,\p\sb{t}
\left(
\int
K\sb{\lambda,\sigma}(x,\theta,\bar y + (y-\bar y)t)
a\sb{\eur Q}(y)d\theta\,dy
\right)
\nonumber \\
&=&
\lambda r
\int
\left\{
\int\sb{0}\sp 1dt\,
\frac{y-\bar y}{r}\,
\lambda\sp{-1}\p\sb{y}
K\sb{\lambda,\sigma}(x,\theta,\bar y + (y-\bar y)t)
\right\}
a\sb{\eur Q}(y)d\theta\,dy.
\end{eqnarray}
The expression in the curly brackets can be treated
as an integral kernel of another Fourier integral operator
of the same order $\mu$
associated to $\eub C$,
and therefore
\[
\abs{\euf{F}\sb{\lambda,\sigma} a\sb{\eur Q}(x)}
\le\lambda r
\const\sigma\sp{1/m}.
\]
Let us mention that
in (\ref{extra-lambda-r})
$\abs{\frac{y-\bar y}{r}}\le\const$
and that the increase in the order of the symbol
due to the derivative $\p\sb{y}$ is compensated by $\lambda\sp{-1}$.
When the derivative $\p\sb{y}$
acts on $\beta(\mathcal{D}(x,\theta,y)/\sigma)$
(which is hidden inside $K\sb{\lambda,\sigma}$),
the contribution is bounded by $\const\sigma\sp{-1}$
and is also compensated by $\lambda\sp{-1}$.
The integration in $t$ is irrelevant.

This completes the proof of Lemma \ref{key-lemma-hardy-l-infinity}.





\bigskip

\begin{center}
ACKNOWLEDGMENT
\end{center}

The author would like to thank Michael Taylor 
for calling the author's attention to the subject of caustics
and for his continuous interest, help, and encouragement.

\bibliographystyle{amsalpha}
\bibliography{all,books,berestycki,brenner,cazenave,comech,greenleaf,hormander,melrose,phong,seeger,shatah,sogge,sugimoto,taylor}
\end{document}